\documentclass{amsart}
\usepackage{amscd,amssymb,hyperref}

\theoremstyle{plain}
\newtheorem{prop}[subsection]{Proposition}
\newtheorem{thm}[subsection]{Theorem}
\newtheorem*{theo}{Theorem}
\newtheorem{lem}[subsection]{Lemma}

\newtheorem{definition}[subsection]{Definition}
\theoremstyle{remark}
\newtheorem{rem}[subsection]{Remark}

\theoremstyle{definition}

\numberwithin{equation}{section}

\renewcommand{\b}[1]{\mathbf{#1}}
\newcommand{\smallmath}[1]{\text{{\footnotesize{$#1$}}}}
\newcommand{\set}[1]{\left\{#1\right\}}

\newcommand{\A}{{\mathcal A}}
\newcommand{\Ai}{{\mathcal A}_\infty}

\newcommand{\CC}{{\mathcal C}}
\newcommand{\cE}{{\mathcal E}}
\newcommand{\FF}{{\mathcal F}}
\newcommand{\cG}{{\mathcal G}}
\newcommand{\LL}{{\mathcal L}}

\newcommand{\cT}{{\mathcal T}}

\newcommand{\Z}{{\mathbb Z}}

\newcommand{\C}{{\mathbb C}}
\newcommand{\CP}{{\mathbb{CP}}}

\newcommand{\T}{{({\mathbb C}^*)^n}}
\newcommand{\bS}{{\mathbb S}}

\newcommand{\bl}{{\boldsymbol{\lambda}}}
\newcommand{\ba}{{\boldsymbol{\a}}}
\newcommand{\bb}{{\boldsymbol{\beta}}}

\newcommand{\sfB}{{\sf B}}
\newcommand{\sfD}{{\sf D}}

\newcommand{\sfb}{{\sf b}}

\newcommand{\sfM}{{\sf M}}

\newcommand{\D}{{\Delta}}
\newcommand{\la}{{\lambda }}
\newcommand{\bul}{{\bullet }}
\newcommand{\p}{\partial}

\renewcommand{\a}{{\alpha }}
\renewcommand{\c}{{\gamma }}
\renewcommand{\L}{{\Lambda }}
\renewcommand{\ll}{{\ell }}

\DeclareMathOperator{\ind}{ind}
\DeclareMathOperator{\dep}{dep}
\DeclareMathOperator{\rank}{rank}
\DeclareMathOperator{\codim}{codim}

\DeclareMathOperator{\ii}{i}

\DeclareMathOperator{\Mat}{Mat}
\DeclareMathOperator{\Aut}{Aut}

\DeclareMathOperator{\Hom}{Hom}

\DeclareMathOperator{\Dep}{Dep}

\begin{document}

\title[Gauss-Manin connections for arrangements, III]
{Gauss-Manin connections for arrangements, III \\ Formal
Connections}
\author[D.~Cohen]{Daniel C.~Cohen$^\dag$}
\address{Department of Mathematics, Louisiana State University,
Baton Rouge, LA 70803}
\email{\href{mailto:cohen@math.lsu.edu}{cohen@math.lsu.edu}}
\urladdr{\href{http://www.math.lsu.edu/~cohen/}
{http://www.math.lsu.edu/\~{}cohen}}
\thanks{{$^\dag$}Partially supported by Louisiana Board of Regents
grant LEQSF(1999-2002)-RD-A-01 and by National Security Agency grant
MDA904-00-1-0038}

\author[P.~Orlik]{Peter Orlik$^\ddag$}
\address{Department of Mathematics, University of Wisconsin,
Madison, WI 53706}
\email{\href{mailto:orlik@math.wisc.edu}{orlik@math.wisc.edu}}
\thanks{{$^\ddag$}Partially supported by National Security Agency
grant MDA904-02-1-0019}

\subjclass[2000]{32S22, 14D05, 52C35, 55N25}

\keywords{hyperplane arrangement, local system, Gauss-Manin
connection}

\begin{abstract}
We study the Gauss-Manin connection for the moduli space of an
arrangement of complex hyperplanes in the cohomology of a complex rank
one local system.  We define formal Gauss-Manin connection matrices in
the Aomoto complex and prove that, for all arrangements and all local
systems, these formal connection matrices specialize to Gauss-Manin
connection matrices.
\end{abstract}


\maketitle

\section{Introduction} \label{sec:intro}
Let $\A=\{H_1,\dots,H_n\}$ be an arrangement of $n$ ordered
hyperplanes in $\C^\ll$, and let $\LL$ be a local system of
coefficients on $\sfM=M(\A)=\C^\ll\setminus\bigcup_{j=1}^n H_j$, the
complement of $\A$.  The need to calculate the local system cohomology
$H^\bul(\sfM;\LL)$ arises in various contexts.  For instance, local
systems may be used to study the Milnor fiber of the non-isolated
hypersurface singularity at the origin obtained by coning the
arrangement, see \cite{CS1, CO2}.  In mathematical physics, local
systems on complements of arrangements arise in the Aomoto-Gelfand
theory of multivariable hypergeometric integrals \cite{AK,Gel1,OT2}
and the representation theory of Lie algebras and quantum groups.
These considerations lead to solutions of the Knizhnik-Zamolodchikov
differential equation from conformal field theory, see \cite{SV,Va}.
Here a central problem is the determination of the Gauss-Manin
connection on $H^\bul(\sfM;\LL)$ for certain discriminantal
arrangements, and certain local systems.

A complex rank one local system on $\sfM$ is determined by a
collection of weights $\bl=(\la_1,\dots,\la_n)\in\C^n$.  Associated to
$\bl$, we have a representation $\rho:\pi_1(\sfM)\to\C^*$, given by
$\c_j\mapsto \exp(-2\pi\ii\la_j)$ for any meridian loop $\c_j$ about
the hyperplane $H_j$ of $\A$, and an associated rank one local system
$\LL$ on $\sfM$.  Parallel translation of fibers over curves in the
moduli space $\sfB$ of all arrangements combinatorially equivalent to
$\A$ gives rise to a flat, Gauss-Manin connection on the vector bundle
over $\sfB$ with fiber $H^q(\sfM;\LL)$.  A number of authors have
considered these Gauss-Manin connections, including Aomoto \cite{Ao},
Schechtman and Varchenko \cite{SV,Va}, Kaneko~\cite{JK}, and Kanarek
\cite{HK}.

For local systems which are nonresonant in the sense of Schechtman,
Terao, and Varchenko \cite{STV}, the Gauss-Manin connection matrices
for general position arrangements were found by Aomoto and Kita
\cite{AK}, and Terao \cite{T1} computed these connection matrices for
a larger class of arrangements.  For such local systems, we determined
these connection matrices for all arrangements in \cite{CO4}.  The aim
of this paper is to remove the nonresonance condition.  We construct
formal Gauss-Manin connection matrices in the Aomoto complex for all
arrangements and prove that these formal connection matrices
specialize to Gauss-Manin connection matrices.

In Section \ref{sec:complexes}, we review the stratified Morse theory
construction \cite{C1,CO1} of a finite cochain complex
$(K^\bul(\A),\D^\bul)$, the cohomology of which is naturally
isomorphic to $H^\bul(\sfM;\LL)$.  This leads to the construction of
the universal complex $(K^\bul_\L,\D^\bul(\b{x}))$ for local system
cohomology, where $\L=\C[x_1^{\pm 1},\dots,x_n^{\pm 1}]$.  The
Orlik-Solomon algebra $A^{\bul}$ is a finite dimensional complex
vector space which models the cohomology of the complement with
trivial local system, see \cite{OT1,OT2}.  The ordering of the
hyperplanes of $\A$ provides the \textbf{nbc} basis for $A^\bul$ over
$\C$.  Weights $\bl$ yield an element $a_{\bl}=\sum \la_j a_j$ in
$A^1$, and multiplication by $a_{\bl}$ gives $A^\bul$ the structure of
a cochain complex.  The resulting cohomology $H^\bul(A^\bul,a_\bl)$ is
a combinatorial analog of $H^\bul(\sfM;\LL)$.  The Aomoto complex
$(A^{\bul}_{R},a_{\b{y}})$ is the universal complex for Orlik-Solomon
cohomology.  Here, $A^{\bul}_{R}$ is a free $R$-module over
$A^{\bul}$, where $R=\C[y_1,\dots,y_n]$.

In Section \ref{sec:moduli}, we recall the moduli space of
arrangements with a fixed combinatorial type $\cT$, following Terao
\cite{T1}.  It is determined by the sets $\ind(\cT)$ and $\dep(\cT)$
of independent and dependent collections of $\ll+1$ element subsets of
hyperplanes in $\A_\infty$, the projective closure of $\A$ in
$\CP^\ell$.  Let $\sfB(\cT)$ be a smooth, connected component of this
moduli space.  There is a fiber bundle $\pi:\sfM(\cT) \to \sfB(\cT)$
whose fibers, $\pi^{-1}(\sfb)=\sfM_{\sfb}$, are complements of
arrangements $\A_{\sfb}$ of type $\cT$.  Since $\sfB(\cT)$ is
connected, $\sfM_{\sfb}$ is diffeomorphic to $\sfM$.

The fiber bundle $\pi:\sfM(\cT) \to \sfB(\cT)$ is locally trivial.
Consequently, given a local system on the fiber, there is an
associated vector bundle $\pi^{q}:\b{H}^{q}\to\sfB(\cT)$, with fiber
$(\pi^{q})^{-1}(\sfb)=H^{q}(\sfM_{\sfb};\LL_{\sfb})$ at
$\sfb\in\sfB(\cT)$ for each $q$, $0\le q \le \ell$.  The transition
functions of this vector bundle are locally constant.  Fixing a
basepoint $\sfb \in\sfB(\cT)$, the operation of parallel translation
of fibers over curves in $\sfB(\cT)$ in the vector bundle
$\pi^{q}:\b{H}^{q}\to\sfB(\cT)$ provides a complex representation
$\Psi^{q}_{\cT}:\pi_{1}(\sfB(\cT),\sfb) \longrightarrow
\Aut_\C(H^{q}(\sfM_{\sfb};\LL_{\sfb}))$.  The loops of primary
interest are those linking moduli spaces of codimension one
degenerations of $\cT$.  Such a degeneration is a type $\cT'$ whose
moduli space $\sfB(\cT')$ has codimension one in the closure of
$\sfB(\cT)$.  In this case, we say that $\cT$ covers $\cT'$.  Let
$\c(\cT') \in \pi_{1}(\sfB(\cT),\sfb)$ be such a loop.

The representation $\Psi^q_{\cT}$ gives rise to a Gauss-Manin
connection.  Let $\Omega^q_{\LL}(\cT',\cT)$ be a connection matrix
associated to the loop $\c(\cT')$.  A key idea in this paper is to
extract information about arbitrary arrangements using information
about general position arrangements, of type $\cG$, whose dependent
set is empty.  The moduli space $\sfB(\cG)$ of general position
arrangements is the complement of a divisor $\sfD(\cG)= \bigcup_J
\sfD_J$ in $(\CP^\ll)^n$.  The components $\sfD_J$ of this divisor are
irreducible hypersurfaces indexed by $\ll +1$ element subsets $J$ of
$[n+1]$.  For each such $J$, let $\cG_J$ denote the combinatorial type
of arrangements with $J$ as the only dependent set.

Schechtman, Terao, and Varchenko \cite{STV}, refining work of Esnault,
Schechtman, and Viehweg \cite{ESV}, found conditions on the weights
which insure that the local system cohomology groups vanish except in
the top dimension.  These conditions depend only on the type $\cT$, so
we call weights satisfying them $\cT$-nonresonant.  Falk and Terao
\cite{FT} constructed a basis for the single nonvanishing cohomology
group, called the $\beta$\textbf{nbc} basis.  The matrices
$\Omega_{\LL}(\cG_J,\cG):=\Omega^\ell_{\LL}(\cG_J,\cG)$ were computed
by Aomoto and Kita \cite{AK}.

In \cite{CO4}, we determined Gauss-Manin connection matrices for all
pairs of arrangement types in the nonresonant case.  Let $\cT$ be a
combinatorial type which covers the type $\cT'$.  Let $\bl$ be a
collection of $\cT$--nonresonant weights which defines the local
system $\LL$.  Weights which are $\cT$-nonresonant are necessarily
$\cG$-nonresonant.  Thus in this case, the only nonvanishing
cohomology groups are in the top dimension.  Since these groups depend
only on the combinatorial type, we write
$H^\ell(\cT;\LL)=H^\ell(\sfM_\sfb;\LL_\sfb)$ and its analog for the
general position type.  Define an endomorphism of $H^\ell(\cG;\LL)$ by
\[
\tilde{\Omega}_{\LL}(\cT',\cT)=\sum_{J \in \dep(\cT',\cT)}
m_J(\cT') \cdot \Omega_{\LL}(\cG_J,\cG).
\]
Here, $m_J(\cT')$ is the order of vanishing of the restriction of a
defining polynomial for $\sfD_J$ to $\overline{\sfB}(\cT)$ along
$\sfB(\cT')$, and $\dep(\cT',\cT)=\dep(\cT') \setminus \dep(\cT)$.  In
\cite{CO4}, we showed that there is a commutative diagram
\begin{equation} \label{eq:gen}
\begin{CD}
H^\ll(\cG;\LL)@>>> H^\ll(\cT;\LL)\\
@VV\tilde{\Omega}_{\LL}(\cT',\cT)V
@VV{\Omega}_{\LL}(\cT',\cT)V \\
H^\ll(\cG;\LL) @>>> H^\ll(\cT;\LL)
\end{CD}
\end{equation}
In particular, if $P$ is the matrix of the surjection $H^\ll(\cG;\LL)
\twoheadrightarrow H^\ll(\cT;\LL)$ in the respective
$\beta$\textbf{nbc} bases, then a Gauss-Manin connection matrix
$\Omega_{\LL}(\cT',\cT)$ is determined by the matrix equation
\begin{equation} \label{eq:gpgm}
P \cdot \Omega_{\LL}(\cT',\cT) = \tilde{\Omega}_{\LL}(\cT',\cT)
\cdot P.
\end{equation}

The aim of this paper is to generalize this result to arbitrary
weights.  For weights $\bl$ which define a nontrivial local system
$\LL$ , the cohomology groups $H^q(\cG;\LL)$ vanish for $q<\ell$.
Thus, a statement similar to (\ref{eq:gen}) in dimension $q<\ll$ is
impossible.  We solve this problem by lifting the Gauss-Manin
connection to the Aomoto complex.

The Orlik-Solomon algebra and the Aomoto complex depend only on
the combinatorial type, so we may label them accordingly.  In
Section \ref{sec:formal}, we define for each subset $S$ of
hyperplanes an endomorphism $\tilde{\omega}_S$ of
$A^{\bul}_{R}(\cG)$, the Aomoto complex of the general position
arrangement. We show that these endomorphisms are cochain maps. In
Section \ref{sec:mult}, we generalize the notion of multiplicity to
all subsets of hyperplanes and the notion of dependent sets to
include sets of all sizes and denote these $\Dep(\cT)$.  In
analogy with $\tilde{\Omega}(\cT',\cT)$, we define
\[
\tilde{\omega}(\cT',\cT)=
\sum_{S\in \Dep(\cT',\cT)}m_S(\cT')\cdot\tilde{\omega}_S.
\]
In order to prove that the endomorphism $\tilde{\omega}(\cT',\cT)$
induces a cochain map $\omega(\cT',\cT)$ on $A^{\bul}_{R}(\cT)$, the
Aomoto complex of type $\cT$, we show in Theorem \ref{thm:ideal} that
this endomorphism preserves the subcomplex corresponding to the
Orlik-Solomon ideal of this type.  This provides the analog of
(\ref{eq:gen}) on the level of Aomoto complexes:
\begin{equation*} \label{eq:aomoto}
\begin{CD}
(A^{\bul}_{R}(\cG),a_{\b{y}}) @>>> (A^{\bul}_{R}(\cT),a_{\b{y}}) \\
@VV\tilde{\omega}(\cT',\cT)V      @VV{\omega}(\cT',\cT)V \\
(A^{\bul}_{R}(\cG),a_{\b{y}}) @>>>
(A^{\bul}_{R}(\cT),a_{\b{y}})
\end{CD}
\end{equation*}
The horizontal maps are explicit surjections provided by the
respective \textbf{nbc} bases.  Given weights $\bl$, the
specialization $\b{y}\mapsto \bl$ in the chain endomorphism
$\omega(\cT',\cT)$ defines endomorphisms
$\omega^q_{\bl}(\cT',\cT):A^q(\cT)\rightarrow A^q(\cT)$ for $0\leq q
\leq \ll$.  In Section \ref{sec:GM}, we prove the main result of this
paper.
\begin{theo}
Let $\sfM$ be the complement of an arrangement of type $\cT$ and let
$\LL$ be~the local system on $\sfM$ defined by weights $\bl$.  Suppose
$\cT$ covers $\cT'$.  Let $\phi^q:K^q\twoheadrightarrow H^q(\sfM,\LL)$
be the natural projection.  Then there is an isomorphism
$\tau^q:A^q\rightarrow K^q$ so that a Gauss-Manin connection
endomorphism $\Omega^q_{\LL}(\cT',\cT)$ in local system cohomology is
determined by the equation
\[
\phi^q\circ\tau^q \circ \omega^q_{\bl}(\cT',\cT) =
\Omega^q_{\LL}(\cT',\cT) \circ \phi^q\circ\tau^q.
\]
\end{theo}

The groups $H^{q}(\sfM;\LL)$ are not known in general. When these
groups can be calculated explicitly, our result yields a matrix
equation analogous to \eqref{eq:gpgm} which determines a
Gauss-Manin connection matrix. We use the \textbf{nbc} basis for
$A^q$ and a suitable choice of basis for $H^q(\sfM;\LL)$, see
Section \ref{sec:ex} for examples.

\section{Cohomology Complexes}
\label{sec:complexes}

For an arbitrary complex local system $\LL$ on the complement of an
arrangement $\A$, we used stratified Morse theory in \cite{C1} to
construct a complex $(K^\bul(\A),\D^\bul)$, the cohomology of which is
naturally isomorphic to $H^\bul(\sfM;\LL)$.  We recall this
construction
in the context of rank one local systems from \cite{C1,CO1}.

Let $\bl=(\la_1,\dots,\la_n)\in\C^n$ be a collection of weights.
Associated to $\la$, we have a rank one representation
$\rho:\pi_1(\sfM) \to \C^*$, given by $\rho(\c_j) =t_j$, where
$\b{t}=(t_1,\dots,t_n)\in (\C^*)^n$ is defined by $t_j=
\exp(-2\pi\ii\la_j)$, and $\c_j$ is any meridian loop about the
hyperplane $H_j$ of $\A$, and a corresponding rank one local system
$\LL=\LL_{\b{t}}=\LL_{\bl}$ on $\sfM$.  Note that weights $\bl$ and
$\bl'$ yield identical representations and local systems if
$\bl-\bl'\in\Z^n$.

We assume throughout that $\A$ contains $\ll$ linearly independent
hyperplanes.  Let $\FF$ be a complete flag (of affine subspaces) in
$\C^\ll$,
\begin{equation*} \label{eq:flag}
\FF:\quad \emptyset = \FF^{-1} \subset \FF^0 \subset \FF^1 \subset
\FF^2 \subset\dots \subset \FF^\ll = \C^\ll,
\end{equation*}
transverse to the stratification determined by $\A$, so that $\dim
\FF^q\cap S_X = q-\codim S_X$ for each stratum, where a negative
dimension indicates that $\FF^q\cap S_X=\emptyset$.  For an explicit
construction of such a flag, see \cite[\S1]{C1}.  Let $\sfM^{q} =
\FF^q \cap \sfM$ for each $q$.  Let $K^q=H^q(\sfM^q,\sfM^{q-1};\LL)$,
and denote by $\D^q$ the boundary homomorphism
$H^{q}(\sfM^{q},\sfM^{q-1};\LL) \to H^{q+1}(\sfM^{q+1},\sfM^{q};\LL)$
of the triple $(\sfM^{q+1},\sfM^q,\sfM^{q-1})$.  The following
compiles several results from \cite{C1}.

\begin{thm}\label{thm:Kdot}
Let $\LL$ be a complex rank one local system on $\sfM$.

\begin{enumerate}
\item \label{item:Kdot1}
For each $q$, $0\le q \le \ll$, we have $H^i(\sfM^q,\sfM^{q-1};\LL)=0$
if $i \neq q$, and $\dim_\C H^q(\sfM^q,\sfM^{q-1};\LL) = b_q(\A)$ is
equal to the $q$-th Betti number of $\sfM$ with trivial local
coefficients $\C$.

\item \label{item:Kdot2}
The system of complex vector spaces and linear maps
$(K^\bul,\D^\bul)$,
\[
K^0 \xrightarrow{\ \D^{0}\ } K^1 \xrightarrow{\ \D^1\ } K^2
\xrightarrow{\phantom{\D^{1}}} \cdots
\xrightarrow{\phantom{\D^{1}}} K^{\ll-1} \xrightarrow{\
\D^{\ll-1}\,} K^\ll,
\]
is a complex $(\D^{q+1}\circ\D^q=0)$.  The cohomology of this complex
is naturally isomorphic to $H^\bul(\sfM;\LL)$, the cohomology of
$\sfM$
with coefficients in $\LL$.
\end{enumerate}
\end{thm}

The dimensions of the terms, $K^{q}$, of the complex
$(K^\bul,\D^\bul)$ are independent of $\b{t}$ (resp., $\bl$, $\LL$).
Write $\D^\bul=\D^\bul(\b{t})$ to indicate the dependence of the
complex on $\b{t}$, and view these boundary maps as functions of
$\b{t}$.  Let $\L=\C[x_1^{\pm 1},\dots,x_n^{\pm 1}]$ be the ring of
complex Laurent polynomials in $n$ commuting variables, and for each
$q$, let $K^q_\L= \L \otimes_{\C} K^q$.

\begin{thm}[{\cite[Thm.~2.9]{CO1}}] \label{thm:univcx}
For an arrangement $\A$ of $n$ hyperplanes with complement $\sfM$,
there exists a universal complex $(K^\bul_\L,\D^\bul(\b{x}))$ with
the following properties:
\begin{enumerate}
\item \label{item:univcx1}
The terms are free $\L$-modules, whose ranks are given by the
Betti numbers of $\sfM$, $K^q_\L \simeq \L^{b_q(\A)}$.

\item \label{item:univcx2}
The boundary maps, $\D^q(\b{x}): K^q_\L \to K^{q+1}_\L$ are
$\L$-linear.

\item \label{item:univcx3}
For each $\b{t}\in\T$, the specialization $\b{x} \mapsto \b{t}$ yields
the complex $(K^\bul,\D^\bul(\b{t}))$, the cohomology of which is
isomorphic to $H^\bul(\sfM;\LL_\b{t})$, the cohomology of $\sfM$ with
coefficients in the local system associated to $\b{t}$.
\end{enumerate}
\end{thm}
The entries of the boundary maps $\D^{q}(\b{x})$ are elements of the
Laurent polynomial ring $\L$, the coordinate ring of the complex
algebraic $n$-torus.  Via the specialization $\b{x} \mapsto \b{t} \in
\T$, we view them as holomorphic functions $\T\to\C$.  Similarly, for
each $q$, we view $\D^{q}(\b{x})$ as a holomorphic map
$\D^{q}:\T\to\Mat(\C)$, $\b{t}\mapsto \D^{q}(\b{t})$.

Recall the Orlik-Solomon algebra of $\A$, $A=A(\A)$.  It is the
quotient of the exterior algebra on generators $a_j$, for $1\leq j\leq
n$ by an ideal defined next.  Given $T=\{j_1,\ldots, j_q\}$, we call
$H_{j_1} \cap \ldots \cap H_{j_q}$ the intersection of $T$.  We call
$T$ {\em dependent} if the intersection of $T$ is not empty and $
\codim(H_{j_1}\cap \ldots \cap H_{j_q}) <q$.  Define $\p :A(\A)
\rightarrow A(\A)$ by $\p(1)=0$, $\p a_j=1$ and for $q>1$
\[
\p a_T=\sum_{k=1}^q(-1)^{k-1}a_{j_1}  \ldots
a_{j_{k-1}} \widehat{a_{j_k}}
a_{j_{k+1}} \ldots a_{j_q}.
\]
The ideal $I(\A)$ is generated by $\p a_T$ where $T$ is dependent
and $a_T$ where the intersection of $T$ is empty.

\begin{rem}\label{rem:vsiso}
The canonical graded algebra isomorphism $H^{\bul}(\sfM,\C)\simeq
A^{\bul}(\A)$ induces an isomorphism of vector spaces $K^q (\A)\simeq
A^q(\A)$ for $0\leq q \leq \ll$.  This isomorphism is not canonical.
\end{rem}

Let $a_{\bl}=\sum_{j=1}^n \la_j a_j \in A^1$.  Since
$a_{\bl}a_{\bl}=0$, $(A^{\bul},a_{\bl})$ is a complex.  There is a
universal complex for the cohomology, $H^{\bul}(A^{\bul},a_{\bl})$.
Let $R=\C[y_{1},\dots,y_{n}]$ be the polynomial ring.  The {\em Aomoto
complex} $(A^{\bul}_{R},a_{\b{y}})$ has terms $A^{q}_{R}=R\otimes_{\C}
A^{q}\simeq R^{b_q(\A)}$, and boundary maps given by
$p(\b{y})\otimes\eta \mapsto \sum y_{j}p(\b{y}) \otimes a_{j} \wedge
\eta$.  For $\bl \in \C^{n}$, the specialization $\b{y}\mapsto \bl$ of
the Aomoto complex $(A^{\bul}_{R},a_{\b{y}})$ yields the Orlik-Solomon
complex $(A^{\bul},a_{\bl})$.  The following result was established in
\cite{CO1}.

\begin{thm} \label{thm:approx}
For any arrangement $\A$, the Aomoto complex
$(A^{\bul}_{R},a_{\b{y}})$ is chain equivalent to the linearization of
the universal complex $(K^{\bul}_{\L},\D^{\bul}(\b{x}))$.
\end{thm}

\section{Moduli Spaces}
\label{sec:moduli}

Let $\cT$ be the combinatorial type of the arrangement $\A$ of $n$
hyperplanes in $\C^\ell$ with $n \ge \ell \geq 1$.  We consider
the family of all arrangements of type $\cT$.  Recall  that $\A$
is ordered by the subscripts of its hyperplanes and we assume that
$\A$, and hence every arrangement of type $\cT$, contains $\ell$
linearly independent hyperplanes.

Choose coordinates $\b{u}=(u_1,\dots,u_\ll)$ on $\C^\ll$.  The
hyperplanes of an arrangement of type $\cT$ are defined by linear
polynomials $f_{i} = b_{i,0} + \sum_{j=1}^{\ell} b_{i,j} u_{j} \,\,(i
= 1,\dots, n)$.  We embed the arrangement in projective space and add
the hyperplane at infinity as last in the ordering, $H_{n+1}$.  The
moduli space of all arrangements of type $\cT$ may be viewed as the
set of matrices
\begin{equation}\label{eq:point}
\sfb=
\begin{pmatrix}
b_{1,0} & b_{1,1} & \cdots & b_{1,\ell}\\
b_{2,0} & b_{2,1} & \cdots & b_{2,\ell}\\
\vdots  & \vdots  & \ddots & \vdots \\
b_{n,0} & b_{n,1} & \cdots & b_{n,\ell}\\
1       & 0       & \cdots & 0
\end{pmatrix}
\end{equation}
whose rows are elements of $\CP^\ell$, and whose
$(\ell+1)\times(\ell+1)$ minors satisfy certain dependency
conditions, see \cite[Prop.~ 9.2.2]{OT2}.

Given a subset $I=\{i_1,\dots,i_{\ell+1}\}$ of
$[n+1]:=\{1,\dots,n,n+1\}$, let $\Delta_{I}=\Delta_{I}(\sfb)$ denote
the determinant of the submatrix of $\sfb$ with rows specified by $I$.
For any combinatorial type $\cT$, let $\ind(\cT)$ denote the set of
all $\ell+1$ element subsets $I$ of $[n+1]$ for which $\Delta_{I}\neq
0$ in type $\cT$.  If $\cT$ is realizable, $\ind(\cT)$ is the set of
all subsets $I$ for which $\{H_{i_1},\dots,H_{i_{\ell+1}}\}$ is
linearly independent in the projective closure of an arrangement $\A$
of type $\cT$.  Similarly, let $\dep(\cT)$ be the set of all $\ell+1$
element subsets $J$ of $[n+1]$ for which $\Delta_{J}=0$ in type $\cT$.
The moduli space of type $\cT$ is
\[
\{\sfb \in (\CP^{\ell})^{n} \mid \Delta_{I}(\sfb) \neq 0 \text{
for } I \in \ind(\cT),\Delta_{J}(\sfb)=0 \text{ for } J \in
\dep(\cT)\}.
\]
Let $\sfB(\cT)$ be a smooth, connected component of this moduli
space. Corresponding to each $\sfb \in \sfB(\cT)$, we have an
arrangement $\A_\sfb$, combinatorially equivalent to $\A$, with
hyperplanes defined by the first $n$ rows of the matrix equation
$\sfb \cdot \tilde{\sf u}=0$, where $\tilde{\sf u} =
\begin{pmatrix} 1 & u_1 & \cdots & u_\ell\end{pmatrix}^\top$.
Let $\sfM_\sfb=M(\A_\sfb)$ be the complement of $\A_\sfb$. Let
\[
\sfM(\cT) = \{ (\sfb,{\sf u}) \in (\CP^\ell)^n \times \C^\ell \mid
\sfb \in \sfB(\cT) \ \hbox{and}\ {\sf u} \in \sfM_\sfb\},
\]
and define $\pi_{\cT}:\sfM(\cT) \to \sfB(\cT)$ by
$\pi_{\cT}(\sfb,{\sf u})=\sfb$.  A result of Randell \cite{Ra}
implies that $\pi_{\cT}:\sfM(\cT) \to \sfB(\cT)$ is a fiber
bundle, with fiber $\pi_{\cT}^{-1}(\sfb) = \sfM_\sfb$.

For each $\sfb \in\sfB(\cT)$, weights $\bl$ define a local system
$\LL_\sfb$ on $\sfM_\sfb$.  Since $\pi_{\cT}:\sfM(\cT) \to
\sfB(\cT)$ is locally trivial, there is an associated vector
bundle $\pi^{q}:\b{H}^{q}\to\sfB(\cT)$, with fiber
$(\pi^{q})^{-1}(\sfb)=H^{q}(\sfM_{\sfb};\LL_{\sfb})$ at
$\sfb\in\sfB(\cT)$ for each $q$, $0\le q \le \ell$. The transition
functions of this vector bundle are locally constant.  Fixing a
basepoint $\sfb\in\sfB(\cT)$, the operation of parallel
translation of fibers over curves in $\sfB(\cT)$ in the vector
bundle $\pi^{q}:\b{H}^{q}\to\sfB(\cT)$ provides a complex
representation
\begin{equation} \label{eq:Hqrep}
\Psi_{\cT}^{q}:\pi_{1}(\sfB(\cT),\sfb) \longrightarrow
\Aut_\C(H^{q}(\sfM_{\sfb};\LL_{\sfb})).
\end{equation}

By Theorem \ref{thm:Kdot}, the local system cohomology of $\sfM_\sfb$
may be computed using the Morse theoretic complex $K^\bul(\A_\sfb)$.
The fundmental group of $\sfB(\cT)$ acts by chain automorphisms on
this complex, see \cite[Cor.~3.2]{CO3}, yielding a representation
\begin{equation*} \label{eq:KArep}
\psi_{\cT}^{\bul}:\pi_{1}(\sfB(\cT),\sfb) \longrightarrow
\Aut_\C(K^{\bul}(\A_{\sfb})).
\end{equation*}

\begin{thm}[\cite{CO3}] \label{thm:inducedrep}
The representation $\Psi^{q}_\cT:\pi_{1}(\sfB(\cT),\sfb) \to
\Aut_\C(H^{q}(\sfM_{\sfb};\LL_{\sfb}))$ is induced by the
representation $\psi^{\bul}_\cT:\pi_{1}(\sfB(\cT),\sfb) \to
\Aut_\C(K^{\bul}(\A_{\sfb}))$.
\end{thm}

The vector bundle $\pi^q:\b{H}^q\to\sfB(\cT)$ supports a Gauss-Manin
connection corresponding to the representation (\ref{eq:Hqrep}).  Over
a manifold $X$, there is a well known equivalence between local
systems and complex vector bundles equipped with flat connections, see
\cite{De,Ko}.  Let $\b{V}\to X$ be such a bundle, with connection
$\nabla$.  The latter is a $\C$-linear map $\nabla:\cE^0(\b{V}) \to
\cE^1(\b{V})$, where $\cE^p(\b{V})$ denotes the complex $p$-forms on
$X$ with values in $\b{V}$, which satisfies $\nabla(f\sigma)= \sigma
df + f \nabla(\sigma)$ for a function $f$ and $\sigma\in\cE^0(\b{V})$.
The connection extends to a map $\nabla:\cE^p(\b{V}) \to
\cE^{p+1}(\b{V})$ for $p\ge 0$, and is flat if the curvature
$\nabla\circ\nabla$ vanishes.  Call two connections $\nabla$ and
$\nabla'$ on $\b{V}$ isomorphic if $\nabla'$ is obtained from $\nabla$
by a gauge transformation, $\nabla'=g\circ\nabla\circ g^{-1}$ for some
$g:X\to\Hom(\b{V},\b{V})$.

The aforementioned equivalence is given by $(\b{V},\nabla) \mapsto
\b{V}^{\nabla}$, where $\b{V}^{\nabla}$ is the local system, or
locally constant sheaf, of horizontal sections $\{\sigma \in
\cE^0(\b{V})\mid \nabla(\sigma)=0\}$.  There is also a well known
equivalence between local systems on $X$ and finite dimensional
complex representations of the fundamental group of $X$.  Note that
isomorphic connections give rise to the same representation.  Under
these equivalences, the local system on $X=\sfB(\cT)$ induced by the
representation $\Psi^q_\cT$ corresponds to a flat connection on the
vector bundle $\pi^q:\b{H}^q\to\sfB(\cT)$, the Gauss-Manin connection.

Let $\c\in\pi_1(\sfB(\cT),\sfb)$, and let $g:\bS^1 \to \sfB(\cT)$ be a
representative loop.  Pulling back the bundle $\pi^q:\b{H}^q \to
\sfB(\cT)$ and the Gauss-Manin connection $\nabla$, we obtain a flat
connection $g^*(\nabla)$ on the vector bundle over the circle
corresponding to the representation of
$\pi_1(\bS^1,1)=\langle\zeta\rangle=\Z$ given by $\zeta \mapsto
\Psi^q_\cT(\c)$.  This vector bundle is trivial since any map from the
circle to the relevant classifying space is null-homotopic.
Specifying the flat connection $g^*(\nabla)$ amounts to choosing a
logarithm of $\Psi^q_\cT(\c)$.  The connection $g^*(\nabla)$ is
determined by a connection $1$-form $dz/z \otimes \Omega^q_\cT(\c)$,
where the connection matrix $\Omega^q_\cT(\c)$ corresponding to $\c$
satisfies $\Psi^q_\cT(\c) = \exp(-2 \pi\ii \Omega^q_\cT(\c))$.  If
$\c$ and $\hat\c$ are conjugate in $\pi_1(\sfB(\cT),\sfb)$, then the
resulting connection matrices are conjugate, and the corresponding
connections on the trivial vector bundle over the circle are
isomorphic.  In this sense, the connection matrix $\Omega^q_\cT(\c)$
is determined by the homology class $[\c]$ of $\c$.

The loops of primary interest are those linking moduli spaces of
codimension one degenerations of $\cT$.  These are types $\cT'$ whose
moduli spaces $\sfB(\cT')$ have codimension one in the closure of
$\sfB(\cT)$.  Define a partial order on combinatorial types as
follows: $\cT \ge \cT' \iff \dep(\cT) \subseteq \dep(\cT')$.  The
combinatorial type $\cG$ of general position arrangements is the
maximal element with respect to this partial order.  Write $\cT >
\cT'$ if $\dep(\cT) \subsetneq \dep(\cT')$.  In this case we define
the relative dependence set $\dep(\cT',\cT)= \dep(\cT') \setminus
\dep(\cT)$.  If $\cT > \cT'$, we say that $\cT$ {\em covers} $\cT'$
and $\cT'$ is a {\em codimension one degeneration} of $\cT$ if there
is no realizable combinatorial type $\cT''$ with $\cT>\cT''>\cT'$.

\begin{lem} \cite[Lemma 5.3]{CO4} \label{lem:cover}
The moduli space $\sfB(\cT')$ has complex codimension one in the
closure $\overline{\sfB}(\cT)$ of the moduli space $\sfB(\cT)$ if
and only if $\cT$ covers $\cT'$.
\end{lem}
Suppose $\cT$ covers $\cT'$, and let $\c \in \pi_1(\sfB(\cT),\sfb)$ be
a simple loop in $\sfB(\cT)$ around a generic point in $\sfB(\cT')$.
Denote a corresponding Gauss-Manin connection endomorphism for
cohomology in the local system $\LL$ by $\Omega_{\LL}(\cT',\cT)$.
These endomorphisms are closely related to certain combinatorial
analogs defined on the Aomoto complex.

\section{Formal Connections}
\label{sec:formal}

Let $(A^{\bul}_{R}(\cG),a_{\b{y}}))$ be the Aomoto complex of a
general position arrangement of $n$ hyperplanes in $\C^\ell$. Let
$T=\{i_1,\ldots, i_q\}$.  If order matters, then we call $T$ a
$q$-tuple and write $T=(i_1,\ldots, i_q)$ and $a_T=a_{i_1} \cdots
a_{i_q}$.  Recall that $H_{n+1}$ is the hyperplane at infinity,
considered the largest in the linear order.  In the formulas
below, a set may contain the index $n+1$ but a tuple may not.  If
$T=(i_1,\ldots,i_q)$ is a $q$-tuple, $1 \leq i_k \leq n$, then
$(j,T)=(j,i_1,\ldots,i_q)$ is the $(q+1)$-tuple which adds $j$
with $1 \leq j \leq n$ to $T$ as its first entry and
$T_k=(i_1,\ldots ,\widehat{i_k},\ldots,i_q)$ is the $(q-1)$-tuple
which deletes $i_k$ from $T$.  We write $S \equiv T$ if $S$ and
$T$ are equal sets.

\begin{definition}
Let $S$ be an index set of size $q+1$. Define an $R$-linear
endomorphism of the Aomoto complex,
$\tilde{\omega}_S:(A^{\bul}_{R}(\cG),a_{\b{y}}) \rightarrow
(A^{\bul}_{R}(\cG),a_{\b{y}})$, as follows. In the formulas below,
$T$ is a $p$-tuple, and $1 \le j \le n$.

If $n+1 \notin S$,
\[
\tilde{\omega}_S^p(a_T)=
\begin{cases}
y_j \p a_{(j,T)} & \text{if $p=q$ and $S\equiv (j,T)$,}\\
a_{\b{y}} \p a_T & \text{if $p=q+1$ and $S\equiv T$,}\\
0 & \text{otherwise.}
\end{cases}
\]

If $n+1 \in S$,
\[
\tilde{\omega}_S^p(a_T)=
\begin{cases}
-\bigl(\sum_{j\in[n]-T} y_j\bigr)a_T
&\text{if $p=q$ and $S\equiv T \cup\{n+1\}$,}\\
(-1)^{k-1}y_j a_{(j,T_k)}
&\text{if $p=q$, $S\equiv (j,T_k)\cup\{n+1\}$, and $j \notin T$,}\\
(-1)^ka_{\b{y}} a_{T_k}
&\text{if $p=q+1$ and $S\equiv T_k \cup \{n+1\}$,}\\
0 & \text{otherwise.}
\end{cases}
\]
\end{definition}

\begin{prop} \label{prop:chain map}
For every $S$, the map $\tilde{\omega}_S$ is a cochain
homomorphism of the Aomoto complex
$(A^{\bul}_{R}(\cG),a_{\b{y}}))$.
\end{prop}
\begin{proof}
Let $S$ be an index set of size $q+1$.  Since
$\tilde{\omega}_S^p=0$ for $p \neq q,q+1$, to show that
$\tilde{\omega}_S$ is a cochain map, it suffices to check
commutativity in the three squares indicated below.
\[
\begin{CD}
\longrightarrow A_R^{q-1}(\cG) @>{a_{\b{y}}}>>
A_R^{q}(\cG) @>{a_{\b{y}}}>>
A_R^{q+1}(\cG) @>{a_{\b{y}}}>>
A_R^{q+2}(\cG) \longrightarrow \\
@VV{\tilde{\omega}_S^{q-1}=0}V  @VV{\tilde{\omega}_S^{q}}V
@VV{\tilde{\omega}_S^{q+1}}V
@VV{\tilde{\omega}_S^{q+2}=0}V \\
\longrightarrow A_R^{q-1}(\cG) @>{a_{\b{y}}}>>
A_R^{q}(\cG) @>{a_{\b{y}}}>>
A_R^{q+1}(\cG) @>{a_{\b{y}}}>> A_R^{q+2}(\cG) \longrightarrow
\end{CD}
\]

If $n+1 \not\in S$, then we may assume that $S= \{1,2,\ldots,q+1\}$.
Otherwise $S=\{U,n+1\}$ and we may assume that $U=
\{n-q+1,\ldots,n\}$.  Since $A^{\bul}_{R}(\cG)$ is free on the
generators $a_T$, we may work with these.

In the first square above, we start with $a_T \in A^{q-1}_R$.  Since
$\tilde{\omega}_S^{q-1}=0$, we need to show that
$\tilde{\omega}^q_S(a_{\b{y}} a_T)=0$.

Suppose $n+1 \not\in S$.  Since $T$ must be equivalent to a subset of
$S$, we may assume that $T=(3,4,\ldots,q+1)$.  Then
$\tilde{\omega}^q_S(a_{\b{y}} a_T)=y_1y_2(\p a_{(1,2,T)}+\p
a_{(2,1,T)})=0$.

Suppose $n+1\in S$.  Since $T$ must be equivalent to a subset of
$U$, we may assume $T=(n-q+2, \ldots,n)$.  Then $a_{\b{y}}
a_T=\sum_{j=1}^{n-q+1}y_j\,a_{(j,T)}$.  We get
\[
\begin{aligned}
\tilde{\omega}^q_S(a_{\b{y}}  a_T) & =\sum_{j=1}^{n-q}y_j
\tilde{\omega}^q_S(a_{(j,T)})+ y_{n-q+1}
\tilde{\omega}^q_S(a_{(n-q+1,T)})\\
&=\sum_{j=1}^{n-q}(-1)^2 y_{n-q+1} y_j a_{(n-q+1,j,T)_2}+
\sum_{j=1}^{n-q}(-1)y_j y_{n-q+1}a_{(j,n-q+1,T)_1}=0.
\end{aligned}
\]

In the second square we start with  $a_T\in A^q_R$.

Suppose $n+1 \not\in S$.  Since $T$ must be equivalent to a subset
of $S$, we may assume that $T=(1,2,\ldots,q)$.  Then
$\tilde{\omega}_S^q(a_T)=y_{q+1} \,\p a_{(q+1,T)}$ and $a_{\b{y}}
\tilde{\omega}_S^q(a_T)= y_{q+1} \,a_{\b{y}} \p a_{(q+1,T)}$.  The
only nonzero term in $\tilde{\omega}_S^{q+1}(a_{\b{y}} a_T)$ is
$\tilde{\omega}_S^{q+1}(y_{q+1}a_{(q+1,T)})=y_{q+1} \,a_{\b{y}} \p
a_{(q+1,T)}$, so the assertion holds.

Suppose $n+1\in S$.  If $T \equiv U$, then
$\tilde{\omega}_S^q(a_T)=(-1)^{1}(\sum_{j=1}^{n-q}y_j)\,
a_{(j,T)_1}$ and $a_{\b{y}} \tilde{\omega}_S^q(a_T)=
-(\sum_{j=1}^{n-q}y_j)\,a_{\b{y}} a_T$. Also, $a_{\b{y}} a_T=
\sum_{j=1}^{n-q}y_j a_{(j,T)}$, so
$\tilde{\omega}_S^{q+1}(a_{\b{y}}
a_T)=-(\sum_{j=1}^{n-q}y_j)\,a_{\b{y}} a_T$ as required.  If $T
\not\equiv U$, then we may assume that $T=(n-q,n-q+1,\ldots,n-1)$.
Then $\tilde{\omega}_S^q(a_T)=(-1)^{2}y_n a_{(n,T)_2}$ so
$a_{\b{y}} \tilde{\omega}_S^q(a_T)=y_n a_{\b{y}} a_{(n,T)_2}$. On
the other hand, there is only one term in $a_{\b{y}} a_T$ on which
$\tilde{\omega}_S^{q+1}$ is nonzero, namely $y_n a_n a_T$.  Since
$\tilde{\omega}_S^{q+1}(y_n a_{(n,T)})=y_n a_{\b{y}} a_{(n,T)_2}$,
the assertion holds.

In the third square we start with $a_T\in A^{q+1}_R$.  Since
$\tilde{\omega}_S^{q+2}=0$, it suffices to show that $a_{\b{y}}
\tilde{\omega}_S^{q+1}(a_T)=0$.  For any $S$, this follows from
$a_{\b{y}} a_{\b{y}}=0$.
\end{proof}

\begin{rem} The map $\tilde{\omega}_S^q$ is given by geometric
considerations in \cite{AK,T1,OT2}.  There are many possible lifts
$\tilde{\omega}_S^{q+1}$ which make $\tilde{\omega}_S$ a cochain map.
However, if we require that $\tilde{\omega}^{q+1}_S(a_T)=0$ unless $T$
is related to $S$ as indicated in the definition, then the lift is
unique.
\end{rem}

A collection of weights $\bl=(\la_1,\dots,\la_n)$ is
$\cG$-nonresonant if $\la_j \notin \Z_{\ge 0}$ and $-\sum_{j=1}^n
\la_j\notin \Z_{\ge 0}$.  The $\beta$\textbf{nbc} basis for the
unique nonvanishing local system cohomology group
$H^\ell(\cG;\LL)$ consists of monomials $\eta_K = \la_{k_1}\cdots
\la_{k_\ll} a_K$, where $K=(k_1,\dots,k_\ll)$, $2 \le k_1 < \dots
< k_\ll \le n$.  In this basis, the Gauss-Manin connection is given by
$\sum d\log \sfD_J \otimes \Omega_{\LL}(\cG_J,\cG)$, where the sum
is over all $\ell+1$ element subsets $J$ of $[n+1]$. The connection
matrices $\Omega_{\LL}(\cG_J,\cG)$ were computed in \cite{OT2},
recovering a result of Aomoto and Kita \cite{AK}.

For $J\subset [n+1]$, let $\tilde\omega_J(\bl)$ be the specialization
$\b{y}\mapsto\bl$ of the endomorphism $\tilde\omega_J$.

\begin{prop}
\label{prop:gp}
If $\bl$ is a collection of $\cG$-nonresonant weights, then for each
$\ell+1$ element subset $J$ of $[n+1]$, the Gauss-Manin connection
matrix $\Omega_{\LL}(\cG_J,\cG)$ is induced by the specialization
$\tilde\omega_J(\bl)$ of the endomorphism $\tilde\omega_J$.
\end{prop}
\begin{proof}
For $\cG$-nonresonant weights $\bl$, the projection $A^\ell(\cG)
\twoheadrightarrow H^\ell(\cG;\LL)$ is given by
\[
a_K\mapsto\frac{1}{\la_{k_1}\cdots \la_{k_\ell}} \eta_K
\text{ if $1 \notin K$ and }
a_K\mapsto-\frac{1}{\la_{k_1}\cdots \la_{k_\ell}}
\sum_{i \notin K} \eta_{(i,K')}
\text{ if $K=(1,K')$.}
\]
For an $\ell+1$ element subset $J$ of $[n+1]$, a calculation with this
projection and the endomorphism $\tilde{\omega}_J^{\ll}(\bl)$ yields
the result.
\end{proof}

\section{Degenerations}
\label{sec:mult}
For a combinatorial type $\cT \neq \cG$, the Orlik-Solomon ideal
$I(\cT)$ gives rise to a subcomplex $I^\bul_R(\cT)$ of the Aomoto
complex $A^{\bul}_{R}(\cG)$, with quotient $A^{\bul}_{R}(\cT)$.  In
order to prove that the endomorphisms $\tilde\omega(\cT',\cT)$ on
$A^{\bul}_{R}(\cG)$ induce endomorphisms $\omega(\cT',\cT)$ on
$A^{\bul}_{R}(\cT)$, we must show that they preserve this subcomplex,
$\tilde\omega(\cT',\cT)(I^\bul_R(\cT)) \subset I^\bul_R(\cT)$.  This
fact is established in the next two sections.

\begin{definition}\label{def:mult}
Given $S \subset [n+1]$, let $N_S(\cT)=N_S(\sfb)$ denote the
submatrix of (\ref{eq:point}) with rows specified by $S$. Let
$\rank N_S(\cT)$ be the size of the largest minor with nonzero
determinant. Define the multiplicity of $S$ in $\cT$ by
$$ m_S(\cT)=|S|-\rank N_S(\cT).$$
\end{definition}
This definition agrees with an interpretation of multiplicity
given above for $S$ with $|S|=\ll+1$. For $1<q\leq n+1$, let
\[
\Dep(\cT)_q=\{\{j_1,\ldots,j_q\}\mid \codim(H_{j_1}\cap \ldots
\cap H_{j_q}) <q\}.
\]
The set $\dep(\cT)$ defined above is  $\Dep(\cT)_{\ll+1}$. Let
$\Dep(\cT)=\cup_q \Dep(\cT)_q$. If $\cT>\cT'$ we define
$\Dep(\cT',\cT)= \Dep(\cT')\setminus \Dep(\cT)$. Let
\[
\tilde{\omega}(\cT)=\sum_{S\in
\Dep(\cT)}m_S(\cT)\cdot\tilde{\omega}_S\qquad \text{ and }\qquad
\tilde{\omega}(\cT',\cT)=\sum_{S\in
\Dep(\cT',\cT)}m_S(\cT')\cdot\tilde{\omega}_S.
\]

For the remainder of the paper, we fix $\cT$ and assume that $\cT'$ is
a codimension one degeneration of $\cT$.  Since the Aomoto complex has
dimension $\ell$, only $|S|\leq \ell +1$ can contribute to the maps
$\tilde{\omega}(\cT',\cT)$.  A {\em circuit} is a minimally dependent
set of hyperplanes.  A generating set for $I(\cT)$ is obtained from
the collection of circuits of $\Ai$.

Fix a circuit $T \in \Dep(\cT)$ with $|T|=q+1$.  If $T=\{U,n+1\}$,
then the hyperplanes of $T$ meet at infinity in $\Ai$, so the
hyperplanes of $U$ have empty intersection in $\A$ and $a_{U}$ is
a generator of $I(\cT)$. If $n+1 \not\in T$, then $\p a_T $ is a
generator of $I(\cT)$. Let
\[
r_T=
\begin{cases}
a_U  & \text{if $T=\{U,n+1\},$}\\
\p a_T  & \text{if $n+1 \not\in T$.}
\end{cases}
\]
Note that $r_T \in A^q_R(\cG)$. We refer to $S \in \Dep(\cT',\cT)$
as $T$-\emph{relevant} if $\tilde{\omega}_S(r_T)\neq 0$. For such
$S$, we have $|S|=q$ or $|S|=q+1$. The next observation follows
from the definition.

\begin{lem}
\label{lem:vanish} Let $T$ be a $q$-tuple and let $S$ be any set.
If $|T \cap S| <q-1$, then $\tilde{\omega}_S (a_T)=0$.
\end{lem}

Terao \cite{T1} classified the three codimension one degeneration
types in the moduli space of an arrangement whose only dependent set
is the circuit $T$.  In Type I, $|S\cap T|<q-1$ for all $S \in
\Dep(\cT',\cT)$.  By Lemma \ref{lem:vanish}, these degenerations are
not $T$-relevant.  In the description of the remaining two types we
use the symbol $T_p^m=T_p \cup\{m\}$ for $ 1\le p \le |T|$ and $m
\not\in T$.
\begin{enumerate}
\item[II:] $\{T_p^m \mid m \not\in T\}$ for each fixed
$ p$, $1\le p \le |T|$,
\item[III:] $ \{T_p^m \mid  1 \le p \le |T| \}$
for each fixed $m \not\in T$.
\end{enumerate}
If $q=1$, then Type II does not appear.

\begin{lem} \label{lem:deps}
Let $T$ be a set of cardinality $q+1$.  If $T_i,T_j \in \Dep(\cT)$
for $i \neq j$, then $T_i \in \Dep(\cT)$ for all $i$, $1\le i \le
q+1$.
\end{lem}
\begin{proof}
Without loss, assume that $T_1,T_2 \in \Dep(\cT)$.  Then there are
nonzero vectors $\ba=(0,\a_2,\a_3,\dots,\a_{q+1})$ and
$\bb=(\beta_1,0,\beta_3,\dots,\beta_{q+1})$ which are annihilated by
the rows of the matrix \eqref{eq:point} indexed by $T_1$ and $T_2$
respectively.  If $\a_i=0$ for some $i \neq 1,2$, then $\a$ is
annihilated by the rows corresponding to $T_i$, hence $T_i \in
\Dep(\cT)$.  If $\a_i \neq 0$, then $\a_i\bb-\beta_i\ba$ is a nonzero
vector annihilated by the rows corresponding to $T_i$, hence $T_i \in
\Dep(\cT)$.
\end{proof}

\begin{prop}
\label{prop:m1} Let $T \in\Dep(\cT)_{q+1}$ be a circuit. Then
there is at most one $j \in T$ so that $T_j \in\Dep(\cT',\cT)$.
\end{prop}
\begin{proof}
We may assume that $T=\{1,\dots,q+1\}$.  Then $T_i=T\setminus\{i\}$.
Since $T$ is a circuit, $T_i$ is independent in type $\cT$ for each
$i\in T$.  Since $\cT$ is the type of an arrangement which contains
$\ell$ linearly independent hyperplanes, there exists a set $J \subset
[n]$ of cardinality $\ell-q$ so that $T_1 \cup J \cup \{n+1\}$ is
independent in $\cT$.  We assert that $T_i \cup J \cup \{n+1\}$ is
independent in $\cT$ for all $i\in T$.

Suppose otherwise.  If, for instance, $T_{q+1} \cup J \cup \{n+1\}
\in\Dep(\cT)$, then there are constants $\a_i$, $\beta_j$, $\xi$, not
all zero so that
\begin{equation} \label{eq:lindep}
\sum_{i=1}^q \a_i \sfb_i + \sum_{j\in J} \beta_j \sfb_j + \xi
\sfb_{n+1}=0,
\end{equation}
where $\sfb_k$ denotes the $k$-th row of the matrix
\eqref{eq:point}.  Since $T$ is a circuit, there are constants
$\zeta_k \neq 0$ so that $\sfb_1 = \sum_{k=2}^{q+1} \zeta_k
\sfb_k$. Substituting this expression in \eqref{eq:lindep} yields
a dependence on the set $T_1 \cup J \cup\{n+1\}$, which is a
contradiction.

Let $S=T\cup J \cup \{n+1\}$.  Then $S_i$ is independent in type $\cT$
for each $i$, $1\le i \le q+1$.  Since $\cT'$ is a codimension one
degeneration of $\cT$, if $T_i \in \Dep(\cT')$, then $\sfB(\cT')$ is
locally defined by the vanishing of $\Delta_{S_i}$ in
$\overline{\sfB}(\cT)$.  If $T_j \in \Dep(\cT')$ for $j \neq i$, then
by Lemma \ref{lem:deps}, $T_k \in \Dep(\cT')$ for every $k$.  So, as
above, $\sfB(\cT')$ is locally defined by the vanishing of
$\Delta_{S_k}$ in $\overline{\sfB}(\cT)$ for every $k$, $1\le k\le
q+1$.  We will show that this is a contradiction by exhibiting a point
in $\overline{\sfB}(\cT)$ for which $\Delta_{S_1}$ vanishes but
$\Delta_{S_k}$ does not vanish for $k \neq 1$.

Assume that $J=\{q+2,\dots,\ell+1\}$.  Then $\sfB(\cT)$ contains
points of the form
\[
\sfb(t)=
\begin{pmatrix}
0 & I_q & 0 \\
0 & v & 0 \\
0 & 0 & I_{\ell-q} \\
F(t) & G(t) & H(t)\\
1 & 0 & 0
\end{pmatrix},
\]
where $I_k$ is the $k \times k$ identity matrix, $v=\begin{pmatrix} t
& 1 & \cdots & 1\end{pmatrix}$, and the submatrix $\begin{pmatrix}F(t)
& G(t) & H(t)\end{pmatrix}$ is chosen so that $\sfb(t)$ satisfies the
dependence and independence conditions of type $\cT$ for each nonzero
$t$.  The point $\sfb(0)$ is in $\overline{\sfB}(\cT)$, and for $1\le
k \le q+1$, $\Delta_{S_k}(\sfb(0))$ vanishes only for $k=1$.  Thus, if
$T_i,T_j \in \Dep(\cT')$, there is a realizable type $\cT''$ such that
$\cT > \cT'' > \cT'$, and $\cT'$ is a degeneration of $\cT$ which is
not of codimesion one.
\end{proof}
This result shows that $\Dep(\cT',\cT)_q$ has no $T$-relevant element
unless $T$ has a codimension one degeneration of Type II. In this
case, there is a unique $p$ so that $T_p$ is the only $T$-relevant
element in $\Dep(\cT',\cT)_q$.  It remains to consider $T$-relevant
$S\in \Dep(\cT',\cT)_{q+1}$.

\begin{lem}
\label{lem:m2} Let $T \in\Dep(\cT)_{q+1}$ be a circuit. If all
$T$-relevant $S \in\Dep(\cT',\cT)_{q+1}$ belong to a family of a
single type, then $m_S(\cT')=1$ for each such $S$.
\end{lem}
\begin{proof}
We may assume that $T=\{U,n+1\}$ where $U=\{1,\ldots,q\}$. Suppose
the degeneration is of Type II so $\{U,k\}\in
\Dep(\cT',\cT)_{q+1}$ for some $k\in [n]-U$.  We argue by
contradiction.  If $m_{\{U,k\}}(\cT')=2$, then in type $\cT'$
there are two linearly independent vectors
$\ba=(\a_{1},\ldots,\a_q,\a_k)$ and
$\bb=(\beta_{1},\ldots,\beta_q,\beta_k)$ which are annihilated by
the rows of (\ref{eq:point}) specified by $\{U,k\}$.  If $\a_1=0$,
then $\{U_1,k\}\in \Dep(\cT')$.  Since $\{U,n+1\}$ is a circuit in
$\cT$, and $\{U_1,k\}\not\in \Dep(\cT)$ by Proposition
\ref{prop:m1}, we have $\{U_1,k\}\in \Dep(\cT',\cT)$ and hence
$\{U_1,k,n+1\}\in \Dep(\cT',\cT)$. This contradicts the assumption
that all $T$-relevant sets $S$ belong to a Type II family.  If
$\a_1\neq 0$, then we use it to eliminate $\beta_1$ and find the
same contradiction.

If the degeneration is of Type III, we may assume that
$\{U_1,p,n+1\}\in \Dep(\cT',\cT)$ with $p \in [n]-U$.  Assuming
that $m_{\{U_1,p,n+1\}}(\cT')=2$ leads to a similar argument.  We
consider the coefficient $\a_{n+1}$ and conclude that
$\{U_1,p\}\in \Dep(\cT',\cT)$ and hence $\{U,p\}\in
\Dep(\cT',\cT)$. This contradicts the assumption that all
$T$-relevant sets $S$ belong to a Type III family.
\end{proof}

\begin{lem}
\label{lem:m3} Let $T \in\Dep(\cT)_{q+1}$ be a circuit. Suppose
$T$ gives rise to codimension one degenerations of both Type II
and Type III. Then the Type II family is unique. For each Type III
family there is a unique $p \in [n+1]-T$ so that $\{T_i,p\}$ is in
both families. We call $\{T_i,p\}$ the intersection of these
families. Moreover, $m_{\{T_i,p\}}(\cT')=2$ for each intersection
and $m_S(\cT')=1$ for all other $T$-relevant $S$ in these
families.
\end{lem}
\begin{proof}
We may assume that $T=\{U,n+1\}$ where $U=\{1,\ldots,q\}$.  If $T$
gives rise to a Type II family, then it is of the form $\{\{U,k\}
\mid k\in [n]-U\}$. By Proposition \ref{prop:m1}, there is a
unique $j$ for which $T_j \in \Dep(\cT')_q$.  We may assume that
$j=q+1$ so that $T_{q+1}=U \in \Dep(\cT')$. If $T$ gives rise to
two different families of Type II, then also some $\{U_i,n+1\} \in
\Dep(\cT')$ contradicting Proposition \ref{prop:m1}.

Suppose there is also a Type III family involving $T$.  (There may
be several Type III families involving $T$, but it will be clear
from the proof that we may consider one Type III family at a
time.)  Let $\{U,p\}$ be the intersection of the given Type II
family and this Type III family.  We show first that
$m_{\{U,p\}}(\cT')=2$. Since $U\in \Dep(\cT',\cT)$, it suffices to
prove that row $p$ is a linear combination of the rows specified
by $U$ in (\ref{eq:point}).  Since $\{U,n+1\}\in \Dep(\cT)$, there
is a vector $\ba=(\a_1,\ldots,\a_q,\a_{n+1},0)$ which is
annihilated by the rows $\{U,n+1,p\}$ of (\ref{eq:point}).  Since
$\{U,n+1\}$ is a circuit, all $\a_i\neq 0$.  This dependency holds
also in type $\cT'$.  Since $\{U_1,p,n+1\}\in \Dep(\cT')$, we also
have a vector $\bb=(0,\beta_2,\ldots,\beta_q,\beta_{n+1},\beta_p)$
annihilated by the rows $\{U,n+1,p\}$ of (\ref{eq:point}) in type
$\cT'$.  We claim that $\beta_p\neq 0$, for otherwise we would
have $\{U_1,n+1\}\in \Dep(\cT',\cT)$, contradicting Proposition
\ref{prop:m1}.  The vector $\beta_{n+1}\ba-\a_{n+1}\bb$ provides
the required dependence.  Hence, $m_{\{U,p\}}(\cT') \ge 2$.
Assuming that $m_{\{U,p\}}(\cT') >2$ contradicts Proposition
\ref{prop:m1}.

The fact that the other multiplicities in these families are $1$ is
established as in the proof of Lemma~\ref{lem:m2}.
\end{proof}

\begin{thm}
\label{thm:types} The endomorphism $\tilde{\omega}(\cT',\cT)$
satisfies $\tilde{\omega}(\cT',\cT)(I^\bul_R(\cT))\subset
I^\bul_R(\cT)$ if and only if for each $K \subset [n]$ and each
circuit $T\in \Dep(\cT)$, we have $\sum\tilde{\omega}_S (a_K
r_T)\in I^\bul_R(\cT)$ where the sum is over $T$-relevant $S$ in a
single type of codimension one degeneration involving $T$.
\end{thm}
\begin{proof}
This is clear if $T$ is involved a single type. If more types
appear, then Lemma \ref{lem:m3} shows that each $S$ is the
intersection of at most two types.  Furthermore, all such
intersections have multiplicity 2, so the corresponding $S$ may be
considered individually in their respective types, each time with
multiplicity 1.
\end{proof}

\section{Ideal Invariance}
\label{sec:invar}

\begin{thm}
\label{thm:ideal} If $\cT$ covers $\cT'$, then
$\tilde{\omega}(\cT',\cT)(I^\bul_R(\cT))\subset I^\bul_R(\cT)$.
\end{thm}
\begin{proof}
It follows from Theorem \ref{thm:types} that we may argue on the
different types independently.  It suffices to show that for every
circuit $T\in \Dep(\cT)$, every $k$-tuple $K$, and every degeneration
$\cT'$ of $\cT$, we have $\tilde{\omega}(\cT')(a_K r_T) \in
I^\bul_R(\cT)$.  There are several cases to consider because the
generators of the ideal are defined in terms of affine coordinates,
while the natural action of the symmetric group on the hyperplanes is
in projective space.  Note that $n+1 \not\in K$ and we agree to use
the same symbol for the underlying set.  Similarly, if $L$ is a set
which does not contain $n+1$, then we write $L$ for the corresponding
tuple in the natural order.

The following identity will be useful in several parts.
If $J \subset [n]$, then
\[
(\sum_{m \in J}y_m a_m)\p a_J=(\sum_{m \in J}y_m) a_J.
\]

\subsection*{Case 1} In this case, $T \in \Dep(\cT)$ is a circuit with
$n+1\in T$.

In this instance, $T=\{U,n+1\}$, and we may assume that
$U=\{n-q+1,\ldots, n\}$.  First assume $T \subset S \in \Dep(\cT')$.
Clearly, $\{K,T\}\in \Dep(\cT')$.  Let $L=[n]\setminus (K \cup U)$.
We get
\[
\tilde{\omega}^{k+q}_{\{K,U,n+1\}}(a_K a_U)=
-(\sum_{j\in L}y_j) a_K a_{U}.
\]
For every $j \in K$, $\{K_j,T\}\in \Dep(\cT')$. Here
\[
\tilde{\omega}^{k+q}_{\{K_j,T\}}(a_K a_U)=
(-1)^j a_{\b{y}} a_{K_j} a_{U}.
\]
Similarly, for every $j \in K$ and every $m\in L$, $\{K_j^m,T\}\in
\Dep(\cT')$. Here
\[
\tilde{\omega}^{k+q}_{\{K_j^m,U,n+1\}}(a_K a_U)=
a_m a_{K_j} a_{U}.
\]
In the remaining parts of this case, we may assume that
$T \not\subset S$ for $S\in\Dep(\cT')$.

\subsubsection*{Case 1.1}
If $|S\cap \{K,U\}|<k+q-1$ for all other $S\in \Dep(\cT')$, then we
are done by Lemma \ref{lem:vanish}.  If there exists $S\in \Dep(\cT')$
with $|S\cap \{K,U\}|\geq k+q-1$ and $T \not\subset S$, then
$S=\{K,T_p^m\}$ with $m \in [n+1]\setminus T$.  The classification
implies that $T_p^m$ is in one of the remaining two types and all the
other members of that type must also be in $\Dep(\cT')$.

\subsubsection*{Case 1.2}
Suppose $T_p^m$ belongs to Type II. Then $p$ is fixed.  If $p\neq
q+1$, then we may assume that $p=1$.  Thus $\Dep(\cT')$ contains
$S_m=\{m,K,T_{1}\}$ for all $m\in L$.  Since every $S_m$ contains
$F=\{K,T_{1}\}$, we conclude that $F \in \Dep(\cT')$.  Here
$\tilde{\omega}^{k+q}_F(a_K a_U)=(-1)^{k+1}a_{\b{y}} a_K a_{U_{1}}$
and
$\tilde{\omega}^{k+q}_{S_m}(a_K a_U)=(-1)^{k} y_m a_m a_K a_{U_1}$.
Thus
\[
(\tilde{\omega}^{k+q}_F+ \sum_{m\in L} \tilde{\omega}^{k+q}_{S_m})
(a_K a_U)= -y_{n-q+1}a_K a_U.
\]
If $p=q+1$, then $\Dep(\cT')$ contains $ S_m=\{m,K,U\}$ for all $m\in
L$.  Since every $S_m$ contains $F=\{K,U\}$, we conclude that $F \in
\Dep(\cT')$.  Here $\tilde{\omega}^{k+q}_F(a_K a_U)= a_{\b{y}} \p(a_K
a_U)$ and $\tilde{\omega}^{k+q}_{S_m}(a_K a_U)=y_m \p a_{(m,K,U)}= y_m
a_K a_U - y_m a_m \p(a_K a_U)$.  Thus
\[
(\tilde{\omega}^{k+q}_F+ \sum_{m \in L}
\tilde{\omega}^{k+q}_{S_m}) (a_K a_U)= (\sum_{j\in [n]}y_j) a_K a_U.
\]

\subsubsection*{Case 1.3}
Suppose $T_p^m$ belongs to Type III. Then $\Dep(\cT')$ contains for
some fixed $m \in L$ the sets $S_p=\{m,K,T_p\}$ for all $p$ with $1
\leq p\leq q+1$.  If $p\neq q+1$, then $\tilde{\omega}^{k+q}_{S_p}(a_K
a_U)=(-1)^{k+p+1} y_m a_{(m,K,U)_{k+p+1}}=(-1)^{k+p+1} y_m a_m a_K
a_{U_p}$.  Thus $\sum_{p=1}^{q}\tilde{\omega}^{k+q}_{S_p}(a_K a_U)=
(-1)^k y_m a_m a_K \p a_U$.  If $p=q+1$, then
$\tilde{\omega}^{k+q}_{\{m,K,U\}}(a_K a_U)=y_m \p
a_{(m,K,U)}=y_ma_Ka_U- y_ma_m(\p a_K)a_U+(-1)^{k+1}y_ma_ma_K\p a_U$.
Thus
\[
\sum_{p=1}^{q+1}\tilde{\omega}^{k+q}_{S_p}(a_K a_U)=
y_m(a_K-a_m\p a_K) a_U.
\]
This completes the argument in \textbf{Case 1}.

\subsection*{Case 2} In this case, $T\in \Dep(\cT)$ is a circuit with
$n+1 \not\in T$.

In this instance, we may assume that $T=\{1,\ldots, q+1\}$.  We note
that $a_T=a_1 (\p a_T)\in I(\cT)$ and hence $\p(a_Ka_T)\in I(\cT)$.
First assume $T \subset S \in \Dep(\cT')$.  Clearly, $\{K,T\}\in
\Dep(\cT')$ and we have
\begin{align*}
    \tilde{\omega}^{k+q}_{\{K,T\}}(a_K\p a_T) & =
    \sum_{j\in T}(-1)^{j-1}\tilde{\omega}^{k+q}_{\{K,T\}}(a_K a_{T_j})
    =\sum_{j\in T}(-1)^{j-1}y_j \p a_{(j,K,T_j)}\\
    &=(-1)^k(\sum_{j\in T}y_j)\p (a_K a_T).
\end{align*}
For $j \in K$, $\{K_j,T\}\in \Dep(\cT')$, but
$\tilde{\omega}^{k+q}_{\{K_j,T\}}(a_K \p a_T)=0$. Let
$L=[n]\setminus (K \cup T)$.  For every $j \in K$ and every $m\in
L$, $\{K_j^m,T\}\in \Dep(\cT')$, but only $m=n+1$ gives a nonzero
term:
\[
\tilde{\omega}^{k+q}_{\{K_j,T,n+1\}}(a_K \p a_T)=
(-1)^k(\sum_{s\in T}y_s)a_{K_j}a_T.
\]
In the remaining parts of this case, we may assume that $T
\not\subset S$
for $S \in \Dep(\cT')$.

\subsubsection*{Case 2.1} If $|S\cap \{K,T_j\}|<k+q-1$ for all other
$S\in \Dep(\cT')$ for all $j$, then we are done by Lemma
\ref{lem:vanish}.  If $T \not\subset S$, then $T_p^m \in \Dep(\cT')$
with $m \in [n+1]\setminus T$.

\subsubsection*{Case 2.2} Suppose $T_p^m$ belongs to Type II. Then $p$
is fixed and we may assume that $p=q+1$.  Thus $S_m=\{m,K,T_{q+1}\}
\in \Dep(\cT')$ for all $ m\not\in T$.  Since every $S_m$ contains
$F=\{K,T_{q+1}\}$, we conclude that $F \in \Dep(\cT')$.  We have
\[
\tilde{\omega}^{k+q}_F(a_K\p a_T) =
(-1)^q a_{\b{y}} \p (a_K a_{T_{q+1}}).
\]
Since $T_{q+1}\in \Dep(\cT')$, we also have $\{K_j,T_{q+1},n+1\}\in
\Dep(\cT')$ for all $j \in K$.  Here
$\tilde{\omega}^{k+q}_{\{K_j,T_{q+1},n+1\}}(a_K\p
a_T)=(-1)^{q+j}a_{\b{y}}a_{K_j}a_{T_{q+1}}$.  Thus
\[
\sum_{j \in K}\tilde{\omega}^{k+q}_{\{K_j,T_{q+1},n+1\}}(a_K\p
a_T)= (-1)^{q+1}a_{\b{y}}(\p a_K)a_{T_{q+1}}.
\]
If $m \neq n+1$, then $\tilde{\omega}^{k+q}_{\{m,F\}} (a_K \p
a_T)= (-1)^q y_m \p(a_{m,K,T_{q+1}})$. Note that $m \not\in T$ by
the classification and $m \not\in K$ follows from the expression.
Thus $m \in L$ and we get
\[
\sum_{m \in L}\tilde{\omega}^{k+q}_{\{m,F\}} (a_K \p a_T)=(-1)^q
\sum_{m \in L}y_m \p(a_{m,K,T_{q+1}}).
\]
Consider $m=n+1$. We have $\tilde{\omega}^{k+q}_{\{F,n+1\}} (a_K
\p a_T) =
\sum_{p=1}^{q+1}(-1)^{p-1}\tilde{\omega}^{k+q}_{\{F,n+1\}} (a_K
a_{T_p})$. For $p\neq q+1$, $\tilde{\omega}^{k+q}_{\{F,n+1\}} (a_K
a_{T_p})=(-1)^{p+q}y_p a_K a_{T_{q+1}}$ and
$\tilde{\omega}^{k+q}_{\{F,n+1\}} (a_K a_{T_{q+1}})=-(\sum_{p
\not\in K\cup T_{q+1}}y_p)a_K a_{T_{q+1}}$. Thus
\[
\tilde{\omega}^{k+q}_{\{F,n+1\}} (a_K \p a_{T})
=(-1)^{q+1}(\sum_{j\in [n] \setminus K} y_j )a_K a_{T_{q+1}}.
\]
Similarly, for every $j \in K$ and $s \in L$,
$\{K_j^s,T_{q+1},n+1\}\in \Dep(\cT')$. Here
$\tilde{\omega}^{k+q}_{\{K_j^s,T_{q+1},n+1\}}(a_K a_{T_p})=
(-1)^{q+j+1}y_sa_sa_{K_j}a_{T_{q+1}}$. Thus
\[
\sum_{s \in L}\sum_{j\in
K}\tilde{\omega}^{k+q}_{\{K_j^s,T_{q+1},n+1\}}(a_K a_{T_p})=
(-1)^q(\sum_{s\in L}y_sa_s)(\p a_K)a_{T_{q+1}}.
\]
Summing over all dependent sets in Type II, we must compute
$\xi(a_K\p{a_T})$, where
\[
\xi=\tilde{\omega}^{k+q}_F+\sum_{j \in
K}\tilde{\omega}^{k+q}_{\{K_j,T_{q+1},n+1\}}+ \sum_{m \in
L}\tilde{\omega}^{k+q}_{\{m,F\}}+\tilde{\omega}^{k+q}_{\{F,n+1\}}+
\sum_{s \in L}\sum_{j\in
K}\tilde{\omega}^{k+q}_{\{K_j^s,T_{q+1},n+1\}}.
\]
We get
\begin{align*}
\xi(a_K\p{a_T})&=
(-1)^q a_{\b{y}} \p (a_K a_{T_{q+1}})
+(-1)^{q+1}a_{\b{y}}(\p a_K)a_{T_{q+1}}\\
&\qquad+(-1)^q \sum_{m\in L}y_m \p(a_{m,K,T_{q+1}})
+(-1)^{q+1}\bigl(\sum_{j\in [n] \setminus K} y_j \bigr)a_K
a_{T_{q+1}}\\
&\qquad +(-1)^q\bigl(\sum_{s\in L}y_s a_s)\bigr(\p a_K)a_{T_{q+1}}\\
&=
(-1)^q\bigl(\sum_{j\in K \cup T_{q+1}}y_j + \sum_{m \in L}y_m
-\sum_{j \in [n]\setminus K}y_j \bigr)a_Ka_{T_{q+1}}\\
&\qquad +(-1)^q\bigr(y_{q+1}a_{q+1}+\sum_{m \in L}y_m a_m-
\sum_{m \in L}y_m a_m \bigl) \p (a_Ka_{T_{q+1}})\\
&\qquad + (-1)^q\bigl(-\sum_{j \in [n]\setminus T_{q+1}}y_ja_j +
\sum_{m \in L}y_m a_m \bigr)(\p a_K)a_{T_{q+1}}\\
&=(-1)^{q+1}y_{q+1}[a_Ka_{T_{q+1}}-a_{q+1}\p (a_Ka_{T_{q+1}})+
a_{q+1}(\p a_K)a_{T_{q+1}}]\\
&=-y_{q+1}a_K\p a_T.
\end{align*}

\subsubsection*{Case 2.3} Suppose $T_p^m$ belongs to Type III. Then
$\Dep(\cT')$ contains for some fixed $m\not\in T$ the sets
$S_p=\{m,K,T_p\}$ for all $p \in T$.  If $m\neq n+1$, then
\[
\sum_{p\in T}\tilde{\omega}^{k+q}_{S_p}(a_K \p a_T)  =
\sum_{p\in T}(-1)^{p-1} y_m \p (a_m a_K a_{T_p})
  = y_m(a_K-a_m \p a_K)\p a_T.
\]
For $m=n+1$, we need the formulas
\begin{equation*}
\tilde{\omega}^{k+q}_{\{K,T_p,n+1\}}(a_K a_{T_s}) =
\begin{cases}
(-1)^{p+s-1}y_s a_K a_{T_p} & \text{if $s\neq p$,}\\
-\bigl(\sum_{j\in [n]\setminus \{K,T_p\}}y_j\bigr)a_K a_{T_p} &
\text{if $s=p$.}
\end{cases}
\end{equation*}
We get $\tilde{\omega}^{k+q}_{\{K,T_p,n+1\}}(a_K \p a_T)  =
(-1)^{p} (\sum_{j\in [n]\setminus K} y_j)a_K a_{T_p}$. Thus
\[
\sum_{p\in T}\tilde{\omega}^{k+q}_{\{K,T_p,n+1\}}(a_K \p a_T)
 = -\bigl(\sum_ {j\in [n]\setminus K}y_j\bigr)a_K \p a_{T}.
\]
This completes the argument in \textbf{Case 2}, and hence the proof of
Theorem \ref{thm:ideal}.
\end{proof}

\section{Gauss-Manin Connections}
\label{sec:GM}

Recall the vector bundle $\pi^{q}:\b{H}^{q}\to\sfB(\cT)$, with
fiber $(\pi^{q})^{-1}(\sfb)=H^{q}(\sfM_{\sfb};\LL_{\sfb})$ at
$\sfb\in\sfB(\cT)$ for each $q$, $0\le q \le \ell$ and its
Gauss-Manin connection defined in Section \ref{sec:moduli}. In
this section we investigate this connection and a combinatorial
analog.

Let $\b{A}^q \to \sfB(\cT)$ be the vector bundle over the moduli space
whose fiber at $\sfb$ is $A^q(\A_\sfb)$, the $q$-th graded component
of the Orlik-Solomon algebra of the arrangement $\A_\sfb$.  The
\textbf{nbc} basis provides a global trivialization of this bundle.
Given weights $\bl$, the cohomology of the complex
$(A^{\bul}(\A_\sfb),a_\bl)$ gives rise to an additional vector bundle
$\b{H}^q(A) \to \sfB$ whose fiber at $\sfb$ is the $q$-th cohomology
group of the Orlik-Solomon algebra, $H^q(A^{\bul}(\A_\sfb),a_\bl)$.
Like their topological counterparts, these combintorial vector bundles
admit flat connections, which we call {\em combinatorial Gauss-Manin
connections}, see \cite[\S5]{CO3}.

Theorem \ref{thm:ideal} implies that we have a commutative diagram
of cochain complexes
\[
\begin{CD}
(I^\bul_R(\cT),a_{\b{y}}) @>\iota>> (A^{\bul}_{R}(\cG),a_{\b{y}})
@>p>> (A^{\bul}_{R}(\cT),a_{\b{y}}) \\
@VV\left.\tilde{\omega}(\cT',\cT)\right|_{I^\bul_R(\cT)}V
@VV\tilde{\omega}(\cT',\cT)V
@VV{\omega}(\cT',\cT)V \\
(I^\bul_R(\cT),a_{\b{y}}) @>\iota>> (A^{\bul}_{R}(\cG),a_{\b{y}})
@>p>> (A^{\bul}_{R}(\cT),a_{\b{y}}) \\
\end{CD}
\]
where $\iota:I^\bul_R(\cT) \to A^\bul_R(\cG)$ is the inclusion,
$p:A^\bul_R(\cG) \to A^\bul_R(\cT)$ is the projection provided by the
respective \textbf{nbc} bases, and $\omega(\cT',\cT)$ is the induced
map.  Given weights $\bl$, the specialization $\b{y}\mapsto \bl$ in
the chain endomorphism $\omega(\cT',\cT)$ defines endomorphisms
$\omega^q_{\bl}(\cT',\cT):A^q(\cT)\rightarrow A^q(\cT)$ for $0\leq q
\leq \ll$.

\begin{thm}
\label{thm:GM} Let $\sfM$ be the complement of an arrangement of
type $\cT$ and let $\LL$ be the local system on $\sfM$ defined by
weights $\bl$. Suppose $\cT$ covers $\cT'$.

\begin{enumerate}
\item Let $\rho^q:A^q\twoheadrightarrow H^q(A^\bul(\cT),a_{\bl})$ be
the natural projection. Then a combinatorial Gauss-Manin
connection endomorphism $\Omega^q_{\CC}(\cT',\cT)$ in
Orlik-Solomon  algebra  cohomology is determined by the equation
\[
\rho^q\circ \omega^q_{\bl}(\cT',\cT) = \Omega^q_{\CC}(\cT',\cT)
\circ \rho^q.
\]
\item Let $\phi^q:K^q\twoheadrightarrow H^q(\sfM,\LL)$ be the natural
projection.  Then there is an isomorphism $\tau^q:A^q\rightarrow K^q$
so that a Gauss-Manin connection endomorphism
$\Omega^q_{\LL}(\cT',\cT)$ in local system cohomology is determined by
the equation
\[
\phi^q\circ\tau^q \circ \omega^q_{\bl}(\cT',\cT) =
\Omega^q_{\LL}(\cT',\cT) \circ \phi^q\circ\tau^q.
\]
\end{enumerate}
\end{thm}
\begin{proof}
For $\cT$-nonresonant weights, it is known that
$H^{q}(\sfM_{\sfb};\LL_{\sfb})\simeq H^q(A^{\bul}(\A_\sfb),a_\bl)$ for
all $q$ \cite{STV}.  Moreover, a Gauss-Manin connection matrix
$\Omega^\ell_{\LL}(\cT',\cT)= \Omega^\ell_{\CC}(\cT',\cT)$ in the
$\beta$\textbf{nbc} basis for the unique nonvanishing local system
cohomology group is induced by the endomorphism
$\tilde\Omega_{\LL}(\cT',\cT)$ of the top cohomology of the complement
of a general position arrangement by \cite[Thm.~7.3]{CO4}, see also
\eqref{eq:gpgm}.  In turn, $\tilde\Omega_{\LL}(\cT',\cT)$ is induced
by the specialization $\tilde\omega_\bl(\cT',\cT)$ of the endomorphism
$\tilde\omega(\cT',\cT)$ of the Aomoto complex of a general position
arrangement, see Proposition \ref{prop:gp}.  Hence, it follows from
Theorem \ref{thm:ideal} that $\Omega^\ell_{\LL}(\cT',\cT)$ is induced
by the specialization $\omega_\bl(\cT',\cT)$ of the endomorphism
$\omega(\cT',\cT)$ of the Aomoto complex of type $\cT$.

Thus for $\cT$-nonresonant weights, the endomorphism
$\omega(\cT',\cT)$ induces both the Gauss-Manin connection matrix and
its combinatorial analog in the $\beta$\textbf{nbc} basis.  Since
$\omega(\cT',\cT)$ is a holomorphic map in the variables $\b{y}$ and
the set of $\cT$-nonresonant weights is open and dense,
$\omega(\cT',\cT)$ induces a Gauss-Manin connection endomorphism for
all weights in either cohomology theory.
\end{proof}

\section{An Example}\label{sec:ex}
\begin{figure}[h]
\setlength{\unitlength}{.45pt}
\begin{picture}(300,130)(-200,-110)
\put(-350,-100){\line(0,1){100}} \put(-420,-100){\line(1,1){100}}
\put(-280,-100){\line(-1,1){100}}\put(-420,-60){\line(1,0){140}}
\put(-320,5){3}\put(-280,-55){4}\put(-390,5){1}\put(-355,5){2}
\put(-357,-122){$\mathcal A$}

\put(-150,-100){\line(0,1){100}}\put(-220,-100){\line(1,1){100}}
\put(-80,-100){\line(-1,1){100}}\put(-185,-100){\line(1,1){100}}
\put(-120,5){3}\put(-80,5){4}\put(-190,5){1}\put(-155,5){2}
\put(-157,-122){${\mathcal A}_{1}$}

\put(250,-100){\line(0,1){100}}\put(180,-100){\line(1,1){100}}
\put(320,-100){\line(-1,1){100}}\put(180,-30){\line(1,0){140}}
\put(280,5){3}\put(320,-25){4}\put(210,5){1}\put(245,5){2}
\put(243,-122){${\mathcal A}_{3}$}

\put(50,-100){\line(0,1){100}}\put(-20,-100){\line(1,1){100}}
\put(-20,-60){\line(1,0){140}}
\put(80,5){3}\put(120,-55){4}\put(40,5){12}
\put(43,-122){${\mathcal A}_{2}$}
\end{picture}
\caption{A Codimension One Arrangement and Three Degenerations}
\label{fig:123}
\end{figure}
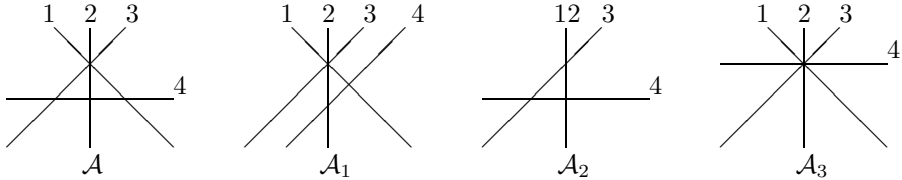

Let $\cT$ be the combinatorial type of the arrangement $\A$ of $4$
lines in $\C^2$ depicted in Figure \ref{fig:123}.  Here $\sfB(\cT)$ is
codimension one in $(\CP^2)^4 =\overline{\sfB}(\cG)$.

For $\cT$-nonresonant weights $\bl$, Gauss-Manin connection
matrices $\Omega^2_\LL(\cT',\cT)$ were determined by Terao
\cite{T1,OT2}. These calculations were recovered in \cite{CO4},
and may also be obtained using the methods of this paper.

For a nontrivial resonant local system $\LL$ on the complement $\sfM$
of $\A$, we used ad hoc methods in \cite[\S6]{CO3} to determine a
Gauss-Manin connection matrix $\Omega^2_\LL(\cT_2,\cT)$, where $\cT_2$
is the combinatorial type of the codimension one degeneration $\A_2$
of $\A$ shown in Figure~\ref{fig:123}.  The results of the previous
sections may be used to obtain these connection matrices for all
codimension one degenerations.  We illustrate this by means of
representative examples.

The universal complex $\bigl(K^\bul_\Lambda,\Delta^\bul(\b{x})\bigr)$
and Aomoto complex $\bigl(A^\bul_R,a_{\b{y}}\bigr)$ of $\A$ are
recorded in \cite[\S6.1]{CO3}.  In particular, the coboundary maps
$\Delta^1(\b{x}):K^1_\L \to K^2_\L$ and $a_{\b{y}}:A^1_R \to A^2_R$
have matrices
\[
\Delta^1 = \Delta^1(\b{x}) =
\smallmath{
\left[\begin{matrix}
x_{3}-x_{2}x_{3} & 1-x_3          & 1-x_{4} & 0       & 0\\
x_1x_{3}-1       & x_{1}-x_{1}x_3 & 0       & 1-x_{4} & 0\\
1-x_{2}          & x_{1}x_2-1     & 0       & 0       & 1-x_{4}\\
0                & 0              & x_1-1   & x_2-1   & x_{3}-1
\end{matrix}\right]}
\]
and
\[
\mu^1 = \mu^1(\b{y}) =
\smallmath{
\left[\begin{matrix}
-y_{2}      & -y_3      & -y_{4} & 0      & 0\\
y_{1}+y_{3} & -y_3      & 0      & -y_{4} & 0\\
-y_{2}      & y_{1}+y_2 & 0      & 0      & -y_{4}\\
0           & 0         & y_1    & y_2    & y_{3}
\end{matrix}\right]}.
\]

The combinatorial types $\cT_i$ of the (multi)-arrangements $\A_i$
shown in Figure \ref{fig:123} are codimension one degenerations of
$\cT$, corresponding to irreducible components of the divisor
$\sfD(\cT)=\overline{\sfB}(\cT) \setminus \sfB(\cT)$.  For these
degenerations, we have
\[
\Dep(\cT_1,\cT)=\{345\},\
\Dep(\cT_2,\cT)=\{12,\ 124,\ 125\},\
\Dep(\cT_3,\cT)=\{124,\ 134,\ 234\}.
\]
The corresponding endomorphisms $\tilde\omega(\cT_i,\cT)$ of the
Aomoto complex $A_R(\cG)$ of a general position arrangement of four
lines are given by
\[
\tilde\omega(\cT_1,\cT)=\tilde\omega_{345},\quad
\tilde\omega(\cT_2,\cT)=\tilde\omega_{12}+\tilde\omega_{124}+
\tilde\omega_{125},\quad
\tilde\omega(\cT_3,\cT)=\tilde\omega_{124}+\tilde\omega_{134}+
\tilde\omega_{234}.
\]

Let $p:A^\bul(\cG) \to A^\bul(\cT)$ be the natural projection, given
in the \textbf{nbc} bases by
\[
p(a_J)=\begin{cases}
a_{1,3} - a_{1,2} &\text{if $J=\{2,3\}$,}\\
a_J &\text{otherwise.}
\end{cases}
\]
By Theorem \ref{thm:ideal}, the chain endomorphisms
$\tilde\omega(\cT_j,\cT)$ of $A^\bul(\cG)$ induce chain endomorphisms
$\omega(\cT_j,\cT)$ of $A^\bul(\cT)$.  A calculation with the
projection $p$ reveals that these are given by
$\omega^1(\cT_j,\cT)=\tilde\omega^1(\cT_j,\cT)$ for each $j$, and
\[
\omega^2(\cT_1,\cT)=
\smallmath{
\left[\begin{matrix}
0 & 0 & 0 & 0 & 0\\
0 & 0 & 0 & 0 & -y_4\\
0 & 0 & 0 & 0 & y_3\\
0 & 0 & 0 & 0 & y_3\\
0 & 0 & 0 & 0 & -y_1-y_2
\end{matrix}\right]}, \quad
\omega^2(\cT_2,\cT)=
\smallmath{
\left[\begin{matrix}
y_1+y_2 & 0 & 0 & 0 & 0\\
y_2 & 0 & 0 & 0 & 0\\
0 & 0 & y_2 & -y_2 & 0\\
0 & 0 & -y_1 & y_1 & 0\\
0 & 0 & 0 & 0 & 0
\end{matrix}\right]},
\]
\[
\omega^2(\cT_3,\cT)=
\smallmath{
\left[\begin{matrix}
y_4 & 0 & -y_4 & y_4 & 0\\
0 & y_4 & -y_4 & 0 & y_4\\
-y_2 & -y_3 & y_2+y_3 & -y_2 & -y_3\\
y_1+y_3 & -y_3 & -y_1 & y_1+y_3 & -y_3\\
-y_2 & y_1+y_2 & -y_1 & -y_2 & y_1+y_2
\end{matrix}\right]}.
\]

\subsection{Nonresonant weights}
For $\cT$-nonresonant weights, local system cohomology is isomorphic
to the cohomology of the Orlik-Solomon complex.  So it suffices to
find combinatorial Gauss-Manin connection matrices in this instance.

Weights $\bl$ are $\cT$-nonresonant if
$\la_1,\la_2,\la_3,\la_4,\la_{\set{1,2,3}},-\la_{[4]} \notin \Z_{\ge
0}$, where $\la_J=\sum_{j\in J}\la_j$ for a set $J$.  The
$\beta$\textbf{nbc} basis for $H^2(A^\bul(\cT),a_\bl)$ is
$\{\eta_{2,4},\eta_{3,4}\}$, where $\eta_{i,j}=\la_i\la_ja_{i,j}$, and
the projection $\rho^2:A^2(\cT) \twoheadrightarrow
H^2(A^\bul(\cT),a_\bl)$ is given by
\[
\rho^2(a_{i,j})=
\begin{cases}
(\la_{\set{1,2}}\eta_{2,4}+\la_2\eta_{3,4})/(\la_1\la_2\la_{\set{1,2,3}})
&\text{if $\{i,j\}=\{1,2\}$,}\\
(\la_3\eta_{2,4}+\la_{\set{1,3}}\eta_{3,4})/(\la_1\la_3\la_{\set{1,2,3}})
&\text{if $\{i,j\}=\{1,3\}$,}\\
-(\eta_{2,4}+\eta_{3,4})/(\la_1\la_4) &\text{if $\{i,j\}=\{1,4\}$,}\\
\ \eta_{i,j}/(\la_i\la_j)&\text{if $\{i,j\}=\{2,4\}$ or $\{3,4\}$.}
\end{cases}
\]

By Theorem \ref{thm:GM}, connection matrices corresponding to the
codimension one degenerations $\cT_j$ of $\cT$ satisfy $\rho^2 \circ
\omega_\bl^2(\cT_j,\cT)= \Omega^2_{\CC}(\cT_j,\cT) \circ \rho^2$.
Calculations with the endomorphisms
$\omega_\bl^2(\cT_j,\cT)=
\left.\omega^2(\cT_j,\cT)\right|_{\b{y}\mapsto\bl}$
and the projection $\rho^2$ yield
\[
\Omega_{\CC}^2(\cT_1,\cT)=
\smallmath{
\left[\begin{matrix}
0 & \la_2 \\ 0 & -\la_{\set{1,2}}
\end{matrix}\right]}, \
\Omega_{\CC}^2(\cT_2,\cT)=
\smallmath{\left[\begin{matrix}
\la_{\set{1,2}} & \la_2 \\ 0 & 0
\end{matrix}\right]}, \
\Omega_{\CC}^2(\cT_3,\cT)=
\smallmath{
\left[\begin{matrix}
\la_{[4]} & 0 \\ 0 & \la_{[4]}
\end{matrix}\right]},
\]
compare \cite[Ex.~10.4.2]{OT2} and \cite[Ex.~8.2]{CO4}.  Similar
calculations recover the remaining connection matrices recorded in
\cite[Ex.~10.4.2]{OT2}.

\subsection{Resonant weights}
A nontrivial resonant local system $\LL$ on $\sfM$ corresponds to a
point $\b{1} \neq \b{t} \in (\C^*)^4$ satisfying $t_1t_2t_3=1$ and
$t_4=1$.  For each such $\b{t}$, we have $H^2(\sfM;\LL)\simeq \C^3$.
Define $\Xi:\L^{5}\to\L^{3}$ and $\Upsilon:R^{5}\to R^{3}$ by
\[
\Xi =
\smallmath{
\left[\begin{matrix}
x_1x_{2}-1 & 0       & x_1x_2+x_4-2       \\
x_{2}-1    & 0       & x_2-1       \\
0          & x_{2}-1 & 0      \\
0          & 1-x_1   & x_3-1 \\
0          & 0       & 1-x_{2}
\end{matrix}\right]}
\  \text{and} \
\Upsilon =
\smallmath{
\left[\begin{matrix}
y_1+y_{2} & 0    & y_1+y_2+y_4    \\
y_{2}     & 0    & y_2    \\
0         & y_2  & 0    \\
0         & -y_1 & y_3  \\
0         & 0    & -y_2
\end{matrix}\right]}.
\]
Note that $\Upsilon$ is the linearization of $\Xi$.  Since $\LL$ is
nontrivial, $t_i \neq 1$ for some $i$.  Assume, without loss, that
$t_2\neq 1$.  For each $\b{t}$ satisfying $t_{1}t_{2}t_{3}=1$,
$t_{4}=1$, and this condition, check that $\rank \Xi(\b{t}) = 3$ and
$\Xi(\b{t}) \circ \D^{1}(\b{t}) = 0$.  So the projection $\C^{5}
\simeq K^{2} \twoheadrightarrow H^{2}(\sfM;\LL) \simeq \C^{3}$ may be
realized as the evaluation $\phi^2=\Xi(\b{t})$.

Via the map $\Upsilon$, the endomorphisms $\omega(\cT_j,\cT)$ induce
maps $\Omega_j:R^3 \to R^3$, which satisfy $\Upsilon \circ
\omega^2(\cT_j,\cT)= \Omega_j \circ \Upsilon$.
Writing $y_J=\sum_{j\in J}y_j$, these maps have matrices
\[
\Omega_1=
\smallmath{
\left[\begin{matrix}
0 & 0 & y_{\set{1,2,4}} \\ 0 & 0 & -y_3 \\ 0 & 0 & -y_{\set{1,2}}
\end{matrix}\right]},
\
\Omega_2=
\smallmath{
\left[\begin{matrix}
y_{\set{1,2}} & 0 & y_{\set{1,2,4}} \\ 0 & y_{\set{1,2}} & -y_3 \\
0 & 0 & 0
\end{matrix}\right]},
\
\Omega_3=
\smallmath{
\left[\begin{matrix}
y_4 & -y_4 & -y_{[4]} \\ -y_{[3]} & y_{[3]} & -y_{[4]} \\
0 & 0 & y_{[4]}
\end{matrix}\right]}.
\]

Let $\bl$ be a collection of weights corresponding to $\b{t}$.  Note
that $\la_2\notin\Z$ since $t_2\neq 1$.  Hence, $\Upsilon(\bl):\C^5
\twoheadrightarrow \C^3$ is surjective for all such $\bl$.  Identify
$H^2(\sfM;\LL) = \C^3$ and let $\tau^2:A^2 \to K^2$ be an isomorphism
for which $\phi^2 \circ \tau^2=\Upsilon(\bl)$.  Theorem \ref{thm:GM}
implies that a Gauss-Manin connection matrix $\Omega^2_\LL(\cT_j,\cT)$
in $H^2(\sfM;\LL)$ corresponding to the degeneration $\cT_j$ of $\cT$
satisfies $\Upsilon(\bl)\cdot \omega^q_{\bl}(\cT',\cT) =
\Omega^q_{\LL}(\cT_j,\cT) \cdot \Upsilon(\bl)$.  Since the equality
$\Upsilon \circ \omega^2(\cT_j,\cT)= \Omega_j \circ \Upsilon$ holds in
the Aomoto complex, the specialization $\b{y}\mapsto\bl$ yields
connection matrices $\Omega_\LL^2(\cT_j,\cT)=\Omega_j(\bl)$.

The endomorphisms $\Omega^1_{\LL}(\cT_j,\cT)$ may be determined by
similar methods. As noted in \cite{CO3}, the endomorphism
$\Omega^1_{\LL}(\cT_2,\cT)$ corresponds to the automorphism of
$H^1(\sfM;\LL) \simeq \C$ given by multiplication by $t_1t_2$, so
has matrix $[\la_{\set{1,2}}]$. (Note that $\la_{\set{1,2}}\in\Z$
if $t_1t_2=1$.)  The endomorphisms $\Omega^1_{\LL}(\cT_1,\cT)$ and
$\Omega^1_{\LL}(\cT_3,\cT)$ are trivial.

Connection matrices $\Omega^q_\LL(\cT',\cT)$ corresponding to other
codimension one degenerations of $\cT$ may be obtained by analogous
calculations.  Note however that the projection $\Upsilon(\bl):A^2 \to
H^2(\sfM;\LL)$ need not be relevant for all degenerations.

\bibliographystyle{amsalpha}

\begin{thebibliography}{00}

\bibitem{Ao} K.~Aomoto, {\em Gauss-Manin connection of integral of
difference products}, J. Math. Soc. Japan \textbf{39} (1987),
191--208;
\href{http://www.ams.org/mathscinet-getitem?mr=88f:32031}
{\MR{88f:32031}}.

\vskip 2pt

\bibitem{AK} K.~Aomoto, M.~Kita, {\em Hypergeometric Functions}
(in Japanese), Springer-Verlag, Tokyo, 1994.

\vskip 2pt

\bibitem{C1} D.~Cohen, {\em Cohomology and intersection cohomology of
complex hyperplane arrangements}, Adv.  Math.  \textbf{97} (1993),
231--266;
\href{http://www.ams.org/mathscinet-getitem?mr=94a:32055}
{\MR{94a:32055}}.

\vskip 2pt

\bibitem{CO1} D. Cohen, P. Orlik, {\em Arrangements and local
systems}, Math. Res. Lett. \textbf{7} (2000), 299--316;
{\href{http://www.ams.org/mathscinet-getitem?mr=2001i:57040}
{\phantom{\ }MR \textbf{2001i}:57040}}.

\vskip 2pt

\bibitem{CO2} \bysame, {\em Some cyclic covers of complements of
arrangements}, Topology Appl.~\textbf{118} (2002), 3--15;
\href{http://www.ams.org/mathscinet-getitem?mr=2003h:32039}
{\MR{2003h:32039}}.

\vskip 2pt

\bibitem{CO3} \bysame, {\em Gauss-Manin connections for
arrangements, {\rm{I}} Eigenvalues}, Compositio Math. \textbf{136} (2003),
299--316.

\vskip 2pt

\bibitem{CO4} \bysame, {\em Gauss-Manin connections for
arrangements, {\rm{II}} Nonresonant weights}, preprint, 2002;
\texttt{\href{http://xxx.lanl.gov/abs/math.AG/0207114}
{math.AG/0207114}}.


\vskip 2pt

\bibitem{CS1} D.~Cohen, A.~Suciu, {\em On Milnor fibrations of
arrangements}, J. London Math. Soc. \textbf{51} (1995), 105--119;
\href{http://www.ams.org/mathscinet-getitem?mr=96e:32034}
{\MR{96e:32034}}.

\vskip 2pt

\bibitem{De} P. Deligne, {\em Equations Diff\'erentielles \`a
Points Singuliers R\'eguliers}, Lect. Notes in Math., vol.~163,
Springer-Verlag, Berlin-New York, 1970.
\href{http://www.ams.org/mathscinet-getitem?mr=54,5232} {\MR{54
\#5232}}

\vskip 2pt

\bibitem{ESV} H.~Esnault, V.~Schechtman, V.~Viehweg, {\em
Cohomology of local systems on the complement of hyperplanes},
Invent. Math. \textbf{109} (1992), 557--561;
\href{http://www.ams.org/mathscinet-getitem?mr=93g:32051}
{\MR{93g:32051}}.  
Erratum, ibid. \textbf{112} (1993), 447;
\href{http://www.ams.org/mathscinet-getitem?mr=94b:32061}
{\MR{94b:32061}}.  

\vskip 2pt

\bibitem{FT} M.~Falk, H.~Terao, {\em $\beta${\bf nbc}-bases for
cohomology of local systems on hyperplane complements},
Trans.~Amer. Math. Soc. \textbf{349} (1997), 189--202.
\href{http://www.ams.org/mathscinet-getitem?mr=97g:52029}
{\MR{97g:52029}}

\vskip 2pt

\bibitem{Gel1}
I. Gelfand, 
{\em General theory of hypergeometric functions}, Soviet Math.
Dokl. \textbf{33} (1986), 573--577;
\href{http://www.ams.org/mathscinet-getitem?mr=87h:22012}
{\MR{87h:22012}}.

\vskip 2pt

\bibitem{GM} M.~Goresky, R.~MacPherson, {\em Stratified Morse
Theory}, Ergeb. Math. Grenzgeb., vol.~14, Springer-Verlag,
Berlin-New~York, 1988;
\href{http://www.ams.org/mathscinet-getitem?mr=90d:57039}
{\MR{90d:57039}}.

\vskip 2pt

\bibitem{HK} H. Kanarek, {\em Gauss-Manin connection arising from
arrangements of hyperplanes}, Illinois J. Math. \textbf{44}
(2000), 741--766;
\href{http://www.ams.org/mathscinet-getitem?mr=2002m:14006}
{\MR{2002m:14006}}.

\vskip 2pt

\bibitem{JK} J. Kaneko, {\em The Gauss-Manin connection of the
integral of the deformed difference product}, Duke Math. J.
\textbf{92} (1998), 355--379;
\href{http://www.ams.org/mathscinet-getitem?mr=99h:32024}
{\MR{99h:32024}}.

\vskip 2pt

\bibitem{Ko} S. Kobayashi, {\em Differential geometry of complex
vector bundles}, Princeton Univ. Press, Princeton, NJ, 1987;
\href{http://www.ams.org/mathscinet-getitem?mr=89e:53100}
{\MR{89e:53100}}.

\vskip 2pt

\bibitem{OT1} P.~Orlik, H.~Terao, {\em Arrangements of
Hyperplanes}, Grundlehren Math. Wiss., vol.~300, Springer-Verlag,
Berlin, 1992;
\href{http://www.ams.org/mathscinet-getitem?mr=94e:52014}
{\MR{94e:52014}}.

\vskip 2pt

\bibitem{OT2} \bysame, {\em Arrangements and Hypergeometric
Integrals}, MSJ Mem., vol.~9, Math.~Soc.~Japan, Tokyo, 2001;
\href{http://www.ams.org/mathscinet-getitem?mr=2003a:32048}
{\MR{2003a:32048}}.

\vskip 2pt

\bibitem{Ra} R. Randell,
{\em Lattice-isotopic arrangements are topologically isomorphic},
Proc. Amer. Math. Soc. \textbf{107} (1989), 555--559;
\href{http://www.ams.org/mathscinet-getitem?mr=90a:57032}
{\MR{90a:57032}}.

\vskip 2pt

\bibitem{STV} V.~Schechtman, H.~Terao, A.~Varchenko, {\em
Cohomology of local systems and the Kac-Kazhdan condition for
singular vectors}, J. Pure Appl. Algebra \textbf{100} (1995),
93--102; \href{http://www.ams.org/mathscinet-getitem?mr=96j:32047}
{\MR{96j:32047}}.

\vskip 2pt

\bibitem{SV} V.~Schechtman and A.~Varchenko, {\em Arrangements of
hyperplanes and Lie algebra homology}, Invent.  Math. \textbf{106}
(1991), 139--194;
\href{http://www.ams.org/mathscinet-getitem?mr=93b:17067}
{\MR{93b:17067}}.

\vskip 2pt

\bibitem{T1} H.~Terao, {\em Moduli space of combinatorially equivalent
arrangements of hyperplanes and logarithmic Gauss-Manin
connections}, Topology Appl.  \textbf{118} (2002), 255--274;
\href{http://www.ams.org/mathscinet-getitem?mr=2003e:32049}
{\MR{2003e:32049}}.

\vskip 2pt

\bibitem{Va} A.~Varchenko, {\em Multidimensional Hypergeometric
Functions and Representation Theory of Lie Algebras and Quantum
Groups}, Adv. Ser. Math. Phys., vol. 21, World Scientific, River
Edge, NJ, 1995;
\href{http://www.ams.org/mathscinet-getitem?mr=99i:32029}
{\MR{99i:32029}}.

\end{thebibliography}

\end{document}